\documentclass[12pt,fullpage]{article}
\usepackage{amsthm}
\usepackage{amsfonts}
\usepackage{xr}
\usepackage{graphicx}
\usepackage{amsmath}
\usepackage{epsfig,subfigure}
\usepackage{color}
\def\R{\mathbb{R}}
\def\N{\mathbb{N}}

\def\epsilon{\varepsilon}

\def\tilde{\widetilde}
\def\div{\mbox{div}}

\newcommand{\SE}{\setcounter{equation}{0} \section}
\newcommand{\be}{\begin{equation}}
\newcommand{\ee}{\end{equation}}
\newcommand{\baa}{\begin{array}}
\newcommand{\eaa}{\end{array}}
\newcommand{\ba}{\begin{eqnarray}}
\newcommand{\ea}{\end{eqnarray}}

\newtheorem{theo}{\bf Theorem}[section]
\newtheorem{lem}[theo]{\bf Lemma}
\newtheorem{pro}[theo]{\bf Proposition}
\newtheorem{cor}[theo]{\bf Corollary}

\newtheorem{rem}[theo]{\bf Remark}

\begin{document}
\date{}
\title{\bf{Large time monotonicity of solutions of reaction-diffusion equations in $\R^N$}}
\author{Emmanuel Grenier$^{\hbox{\small{ a}}}$ and Fran{\c{c}}ois Hamel$^{\hbox{\small{ b}}}$ \thanks{This work has been carried out in the framework of Archim\`ede Labex (ANR-11-LABX-0033) and of the A*MIDEX project (ANR-11-IDEX-0001-02), funded by the ``Investissements d'Avenir" French Government programme managed by the French National Research Agency (ANR). The research leading to these results has also received funding from the ANR within the project NONLOCAL ANR-14-CE25-0013 and from the European Research Council under the European Union's Seventh Framework Programme (FP/2007-2013) / ERC Grant Agreement n.321186 - ReaDi - Reaction-Diffusion Equations, Propagation and Modelling.}\\
\\
\footnotesize{$^{\hbox{a }}$Ecole Normale Sup\'erieure de Lyon, CNRS, UMPA, UMR 5669}\\
\footnotesize{46, all\'ee d'Italie, 69364 Lyon Cedex 07, France}\\
\footnotesize{$^{\hbox{b }}$Aix Marseille Universit\'e, CNRS, Centrale Marseille, I2M, UMR 7373, 13453 Marseille, France}}
\maketitle

\begin{abstract}
In this paper, we consider nonnegative solutions of spatially heterogeneous Fisher-KPP type reaction-diffusion equations in the whole space. Under some assumptions on the initial conditions, including in particular the case of compactly supported initial conditions, we show that, above any arbitrary positive value, the solution is increasing in time at large times. Furthermore, in the one-dimensional case, we prove that, if the equation is homogeneous outside a bounded interval and the reaction is linear around the zero state, then the solution is time-increasing in the whole line at large times. The question of the monotonicity in time is motivated by a medical imagery issue.
\end{abstract}


\SE{Introduction and main results}\label{intro}

In this paper, we consider the Cauchy problem for the following reaction-diffusion equation set in the whole space $\R^N$
\be\label{eq}\left\{\baa{rcll}
u_t & = & \div(A(x)\nabla u)+f(x,u), & t>0,\ x\in\R^N,\vspace{3pt}\\
u(0,x) & = & u_0(x).\eaa\right.
\ee
Here $u_t$ stands for $u_t(t,x)=\frac{\partial u}{\partial t}(t,x)$ and the divergence and the gradient act on the spatial variables $x$. We are interested in the monotonicity in time for large times, when the initial condition is localized and equation~\eqref{eq} is of the monostable Fisher-KPP type. More precisely, the assumptions are listed below.

\subsubsection*{Framework and main assumptions}

The initial condition $u_0$ is in $L^{\infty}(\R^N)$ with $0\le u_0(x)\le 1$ a.e. in $\R^N$ and $u_0$ is non-trivial, in the sense that $\|u_0\|_{L^{\infty}(\R^N)}>0$. We also assume that either there exists $\beta>0$ such that
\be\label{u0gaussian}
u_0(x)=O(e^{-\beta|x|^2})\ \hbox{ as }|x|\to+\infty
\ee
(a particular important case is when $u_0$ is compactly supported), or there exist $0<\gamma\le\delta$ and~$\lambda>0$ such that
\be\label{u0}
\gamma\,e^{-\lambda|x|}\le u_0(x)\le\delta\,e^{-\lambda|x|}\ \hbox{ for all }|x|\hbox{ large enough},
\ee
where $|\cdot|$ denotes the Euclidean norm in $\R^N$.

The diffusion term $A$ is assumed to be a symmetric matrix field $A=(A_{ij})_{1\le i,j\le N}$ of class~$C^{1,\alpha}(\R^N)$ for some $0<\alpha<1$ and uniformly definite positive: there exists a constant~$\nu\ge1$ such that
\be\label{hypA}
\nu^{-1}I\le A(x)\le\nu I\ \hbox{ for all }x\in\R^N,
\ee
in the sense of symmetric matrices, where $I\in\mathbb{S}_N(\R)$ is the identity matrix. One also assumes that~$A$ is locally asymptotically homogeneous at infinity, in the sense that
\be\label{oscA}
\forall\,1\le i,j\le N,\ \ |\nabla A_{ij}(x)|\to0\ \hbox{ as }|x|\to+\infty.
\ee
A particular example of a $C^{1,\alpha}(\R^N)$ matrix field satisfying~\eqref{oscA} is when $A_{ij}(x)$ converges to a constant as $|x|\to+\infty$ for every $1\le i,j\le N$. An important subcase is that of a matrix field~$A$ which is independent of $x$. Notice that, since $A$ is of class $C^{1,\alpha}(\R^N)$, the condition~\eqref{oscA} is equivalent to the fact that the local oscillations of the functions $A_{ij}$ converge to $0$ at infinity, that is, for every $R>0$ and $1\le i,j\le N$,
$$\mathop{\rm{osc}}_{\overline{B(x,R)}}A_{ij}:=\max_{\overline{B(x,R)}}A_{ij}-\min_{\overline{B(x,R)}}A_{ij}\ \to\ 0\ \hbox{ as }|x|\to+\infty,$$
where $B(x,R)$ denotes the open Euclidean ball of center $x$ and radius $R$. However, notice that the matrix fields $A(x)$ satisfying this property may not converge as $|x|\to+\infty$ in general, even in dimension $N=1$.

The reaction term $f:\R^N\times[0,1]\to\R$ is a continuous function, of class $C^{0,\alpha}$ in $x$ uniformly with respect to $u\in[0,1]$, and Lipschitz continuous in $u$, uniformly with respect to $x\in\R^N$. Throughout the paper, one assumes that
\be\label{f1}
f(x,0)=f(x,1)=0\ \hbox{ for every }x\in\R^N
\ee
and that
\be\label{f2}
u\mapsto\frac{f(x,1-u)}{u}\ \hbox{ is non-increasing in }(0,1]
\ee
for every $x\in\R^N$. One also assumes that there exist $\mu>0$ and $s_0\in(0,1)$ such that
\be\label{hypf}
f(x,s)\ge\mu\,s\ \hbox{ for all }(x,s)\in\R^N\times[0,s_0].
\ee
These assumptions imply in particular that $f$ is positive in $\R^N\times(0,1)$ and even that $\inf_{x\in\R^N}f(x,s)\ge\mu s>0$ for every $s\in(0,s_0]$ and $\inf_{x\in\R^N}f(x,s)\ge\mu s_0(1-s)/(1-s_0)>0$ for every $s\in[s_0,1)$. Furthermore, $f$ is assumed to be of class $C^1$ with respect to $u$ in~$\R^N\times([0,s_0]\cup[s_1,1])$ for some $s_1\in(0,1)$ with $f_u=\frac{\partial f}{\partial u}$ bounded and uniformly continuous in $\R^N\times([0,s_0]\cup[s_1,1])$, and of class $C^{0,\alpha}$ with respect to $x$ uniformly in $s\in[0,s_0]\cup[s_1,1]$. Lastly, one assumes that~$f_u(\cdot,0)$ is locally asymptotically homogeneous at infinity, in the sense that, for every $R>0$,
\be\label{oscf}
\mathop{\rm{osc}}_{\overline{B(x,R)}}f_u(\cdot,0)\to0\ \hbox{ as }|x|\to+\infty.
\ee
Notice that~\eqref{oscf} holds if $f_u(\cdot,0)\in C^1(\R^N)$ and $|\nabla f_u(x,0)|\to0$ as $|x|\to+\infty$ or if~$f_u(x,0)$ converges to a constant as~$|x|\to+\infty$ (in particular, if $f_u(\cdot,0)$ is constant). An important class of examples of functions $f$ satisfying the aforementioned hypotheses is when $f(x,u)=r(x)\,g(u)$, where $g$ is of class $C^1$, concave in $[0,1]$, positive in $(0,1)$ with~$g(0)=g(1)=0$, and $r$ is of class~$C^{0,\alpha}(\R^N)$, locally asymptotically homogeneous at infinity and $0<\inf_{\R^N}r\le\sup_{\R^N}r<+\infty$. The archetype is the homogeneous logistic Fisher-KPP~\cite{f,kpp} reaction $f(x,u)=u(1-u)$ with $r(x)=1$ and $g(u)=u(1-u)$ as above. However, for general functions $f(x,u)$ satisfying the above assumptions, slow oscillations at infinity are not excluded, even in dimension $N=1$ (see~\cite{ggn} for the study of one-dimensional equations of the type~\eqref{eq} with slow oscillations as $x\to\pm\infty$).

From the parabolic regularity theory, the solution $u$ of~\eqref{eq} is well-defined for all $t>0$ and it is classical in~$(0,+\infty)\times\R^N$ with
\be\label{0u1}
0<u(t,x)<1\ \hbox{ for all }t>0\hbox{ and }x\in\R^N,
\ee
by the strong parabolic maximum principle. From the assumptions made on $f$, even without~\eqref{oscf}, it is shown in~\cite{bhr} that any stationary solution $p(x)$ of~\eqref{eq} such that $0\le p\le 1$ in~$\R^N$ is either identically equal to $0$ in $\R^N$ or is bounded from below by a positive constant in~$\R^N$. Since $\inf_{x\in\R^N}f(x,s)>0$ for every $s\in(0,1)$, it then follows immediately in the latter case that $p$ is identically equal to $1$ in $\R^N$. Therefore, again from~\cite{bhr}, the solution $u$ of~\eqref{eq} satisfies $u(t,x)\to1$ as $t\to+\infty$ locally uniformly in $x\in\R^N$.

Lastly, from~\cite{bhn}, it is also known that there is $c>0$ such that
\be\label{defc}
\min_{|x|\le c t}u(t,x)\to1\ \hbox{ as }t\to+\infty.
\ee
In other words, the state $1$ invades the whole space as $t\to+\infty$ with at least a positive spreading speed $c>0$. But, the asymptotic spreading speed of $u$ may not be unique, in the sense that some oscillations of the spreading rates of the level sets of $u$ between two different positive speeds are possible in general even for compactly supported initial conditions, see~\cite{ggn}. This means that, in general, there is no speed $c_0>0$ such that~\eqref{defc} holds for all $c\in[0,c_0)$ and $\max_{|x|\ge ct}u(t,x)\to0$ as $t\to+\infty$ for all $c>c_0$. However, when the equation~\eqref{eq} is homogeneous and the initial condition is compactly supported, there exists such a positive spreading speed $c_0$, see e.g.~\cite{aw}.

\subsubsection*{Main results}

The main result of our paper is the following asymptotic time-monotonicity of the solutions of~\eqref{eq}.

\begin{theo}\label{th1}
Under the above assumptions~\eqref{u0gaussian} or~\eqref{u0}, and~\eqref{hypA}-\eqref{oscf}, the solution $u$ of~\eqref{eq} satisfies
\be\label{ut}
\inf_{x\in\R^N}u_t(t,x)\to0\ \hbox{ as }t\to+\infty.
\ee
Furthermore, for every $0<\epsilon<1$, there is a time $T_{\epsilon}>0$ such that
\be\label{Teps}
\forall\,(t,x)\in[T_{\epsilon},+\infty)\times\R^N,\ \ u(t,x)\ge\epsilon\ \Longrightarrow\ u_t(t,x)>0.
\ee
\end{theo}

Property~\eqref{Teps} means the monotonicity in time at large times in the time-dependent sets where $u$ is bounded away from $0$. On the other hand, in the sets where, say, $t\ge 1$ and $u$ is close to $0$, then $u_t$ is close to $0$ too.\footnote{Indeed, if $u(t_n,x_n)\to0$ with $(t_n,x_n)\in[1,+\infty)\times\R^N$, then the functions $v_n(t,x):=u(t+t_n,x+x_n)$ converge locally in $C^{1,2}_{t,x}((-1,+\infty)\times\R^N)$, up to extraction of a subsequence, to a solution $v$ of an equation of the type $v_t=\div(A_{\infty}(x)\nabla u)+f_{\infty}(x,u)$ for some diffusion and reaction coefficients $A_{\infty}$ and $f_{\infty}$ satisfying the same type of assumptions as $A$ and $f$. Furthermore, $v(0,0)=0$ and $0\le v\le 1$ in $(-1,+\infty)\times\R^N$, whence $v=0$ in $(-1,0]\times\R^N$ from the strong maximum principle and then $v=0$ in $(-1,+\infty)\times\R^N$ from the uniqueness of the solutions of the associated Cauchy problem. Finally, $v_t(0,0)=0$ and $u_t(t_n,x_n)=(v_n)_t(0,0)\to v_t(0,0)=0$ as $n\to+\infty$.} Therefore, property~\eqref{Teps} easily yields~\eqref{ut}. Lastly, since
\be\label{conv0}
u(t,x)\to0\hbox{ as }|x|\to+\infty\hbox{ locally uniformly in }t\in[0,+\infty),
\ee
as will be easily seen in the proof of Theorem~\ref{th1} (more precisely, see the proof of Lemma~\ref{lem1} below), property~\eqref{Teps} implies that, for every $T\ge T_{\epsilon}$, the set $\big\{(t,x)\in[T_{\epsilon},T]\times\R^N,\ u(t,x)\ge\epsilon\big\}$ is compact, whence
$$\min_{(t,x)\in[T_{\epsilon},T]\times\R^N,\,u(t,x)\ge\epsilon}u_t(t,x)>0.$$

Let us now comment some earlier related references in the literature. In~\cite{r}, the question of the time-monotonicity at large times had been addressed for the solutions of some reaction-diffusion equations in straight infinite cylinders with advection shear flows and with $f$ being independent of the unbounded variable. Other time-monotonicity results have been obtained in~\cite{bh} for time-global transition fronts of space-heterogeneous reaction-diffusion equations of the type~\eqref{eq} connecting two stable limiting points. In~\cite{z}, the time-monotonicity of the solutions~$u$ of equations $u_t=\Delta u+f(x,u)$ with reactions $f$ of the ignition type or involving a weak Allee effect has been established for large times in the set where $0<\epsilon\le u(t,x)\le 1-\epsilon<1$, for any $\epsilon>0$ small enough. Lastly, we refer to~\cite{dm} for some results on time-monotonicity for small $t$ and large $x$ for the solutions of the homogeneous equation $u_t=\Delta u+g(u)$ which are initially compactly supported.

For the heterogeneous Fisher-KPP type equation~\eqref{eq}, we conjecture that, under the assumptions of Theorem~\ref{th1}, $u_t(t,\cdot)>0$ in $\R^N$ for $t$ large enough. This is still an open question. However, we can answer positively under some additional assumptions on~\eqref{eq} in dimension~$1$.

\begin{theo}\label{th2}
In addition to~\eqref{u0gaussian} or~\eqref{u0},~\eqref{hypA} and~\eqref{f1}-\eqref{hypf}, assume that $N=1$, that~$A'(x)=0$ for $|x|$ large enough and that there are $\lambda^{\pm}>0$, $\theta\in(0,1)$ and two functions $f^{\pm}:[0,1]\to\R$ such that $f(x,u)=f^{\pm}(u)$ for $\pm x$ large enough and $f^{\pm}(u)=\lambda^{\pm}u$ for all $u\in[0,\theta]$. Then there is $\tau>0$ such that the solution $u$ of~\eqref{eq} satisfies
\be\label{utT}
u_t(t,x)>0\ \hbox{ for all }t\ge\tau\hbox{ and }x\in\R.
\ee
\end{theo}

Let us now describe the main ideas of the proof of Theorems~\ref{th1} and~\ref{th2} and the outline of the paper. In Section~\ref{sec2}, the solution $u$ is proved to be $T$-monotone in time ($u(t+T,x)\ge u(t,x)$) at large time $t$ and for all $T$ large enough, by using the decay of $u_0$ at infinity and some Gaussian estimates for the fundamental solution associated with the linear equation obtained from~\eqref{eq}. This $T$-monotonicity is then improved in Section~\ref{sec3} by compactness arguments in the region where $u$ is away from $0$ and from $1$ and then in Section~\ref{sec4} by using in particular the assumption~\eqref{f2} and by an application of the maximum principle in some sets which are defined recursively. In Section~\ref{sec5}, the monotonicity in time is proved in the region where $u$ is close to $1$ by using Harnack inequality applied to the function $1-u$ and some passage to the limit. In Section~\ref{sec6}, the $\tau$-monotonicity in time, for any $\tau>0$, is shown in the region where $u$ is close to $0$, by using some Gaussian estimates as well as some new quantitative inequalities for the fundamental solutions associated with families of linear equations similar to~\eqref{eq} (these new estimates are proved in Section~\ref{secpro1}). Section~\ref{sec7} is devoted to the proof of properties~\eqref{ut} and~\eqref{Teps} of Theorem~\ref{th1}. Lastly, Section~\ref{sec9} is concerned with the proof of Theorem~\ref{th2}, where explicit estimates of the Green function associated to some one-dimensional initial and boundary value problem in half-lines are used.

\begin{rem}{\rm Assume in this remark that, instead of the whole space $\R^N$, equation~\eqref{eq} is set on a smooth bounded domain $\Omega\subset\R^N$ with Neumann type boundary conditions $\mu(x)\cdot\nabla u(t,x)=0$ on $\partial\Omega$, where $\mu$ is a continuous vector field such that $\mu(x)\cdot\nu(x)>0$ for all $x\in\partial\Omega$ and $\nu$ denotes the outward normal vector field on $\partial\Omega$. Then it follows from the arguments used in the proof of Theorem~\ref{th1} (see especially Section~\ref{sec5}) that, under assumptions~\eqref{hypA} and~\eqref{f1}-\eqref{hypf}, any solution $u$ with a nontrivial initial condition $0\le,\not\equiv u_0\le,\not\equiv 1$ is increasing in time in the whole set $\overline{\Omega}$ at large times.}
\end{rem}

\subsubsection*{Modeling and background}

The question of the monotonicity of the solution for large times comes from a simple medical imagery question. A natural way to model a tumor is to introduce a function $\phi(t,x)$ describing the density of tumor cells. In some types of cancers, tumor cells migrate and multiply. They migrate randomly and multiply according to logistic type laws. The simplest model of tumor is therefore the classical KPP equation, as described by Murray~\cite{m}
$$
\phi_t - \nu \Delta \phi = \lambda \phi ( 1 - \phi),
$$
with positive coefficients $\nu$ and $\lambda$. Treatments like radiotherapy or chemotherapy induce the death of a part of tumor cells. A simple way to model a treatment at time $t_0$ is to say that $\phi$ is discontinuous at $t_0$ and
$$
\phi(t_0^+,x) = \beta\,\phi(t_0^-,x) 
$$
for all $x$ and for some $0 < \beta < 1$. 
Now the tumor size can be evaluated through medical imagery devices which detect tumor cells only if their density is large enough, above some threshold $\sigma>0$. The measured size of the tumor is therefore
$$
S(t) = \int_{\R^N} 1_{\phi(t,x) > \sigma} dx .
$$
A natural question is to know whether $S(t)$ can decrease just after a treatment, namely: can the observed size of a tumor decrease whereas its actual total mass $\int_{\R^N}\phi(t,x) dx$ increases~?

Let us detail the link between this question and the positivity of $\phi_t$. For this let~$\Omega(t) = \{ x\in\R^N;\ \phi(t,x) > \sigma \}$, and let $x_0 \in \partial \Omega(t_0^+)$. We have 
$$\baa{rcl}
\phi_t(t_0^+,x_0) & = & \nu \Delta \phi(t_0^+,x_0)+\lambda \phi(t_0^+,x_0) (1 - \phi(t_0^+,x_0))\vspace{3pt}\\
& = & \nu \beta \Delta \phi(t_0^-,x_0)+\lambda \beta \phi(t_0^-,x_0) (1 - \beta \phi(t_0^-,x_0))\vspace{3pt}\\
& = & \beta\phi_t(t_0^-,x_0) - \lambda \beta \phi(t_0^-,x_0) (1 - \phi(t_0^-,x_0)) + \lambda \beta \phi(t_0^-,x_0) (1 - \beta \phi(t_0^-,x_0))\vspace{3pt}\\
& = & \beta \phi_t(t_0^-,x_0) + \lambda \beta (1- \beta) \phi^2(t_0^-,x_0).\eaa$$
The second term is positive, hence if $\phi_t(t_0^-,x) > 0$ everywhere on $\partial\Omega(t_0^+)$, this implies that~$\phi_t(t_0^+,\cdot)$ is positive on $\partial \Omega(t_0^+)$, hence that $S(t)$ is increasing just after $t_0$. The medical imagery question therefore reduces to the study of the sign of $\phi_t$.


\SE{$T$-monotonicity in time}\label{sec2}

Throughout this section and the next ones, one assumes that the conditions~\eqref{hypA}-\eqref{oscf} are fulfilled and $u$ denotes a solution of~\eqref{eq} with initial condition $u_0$ having Gaussian decay at infinity as in~\eqref{u0gaussian} or satisfying~\eqref{u0}. The first step in the proof of Theorem~\ref{th1} consists in showing that~$u$ is $T$-monotone in time.

\begin{lem}\label{lem1}
There is $T>0$ such that
\be\label{defT}
u(1+t,x)\ge u(1,x)\ \hbox{ for all }t\ge T\hbox{ and }x\in\R^N.
\ee
\end{lem}

\noindent{\bf{Proof.}} First of all, as already emphasized, the strong maximum principle implies that~$u(1,x)<1$ for all $x\in\R^N$. Remember also that $u(1,\cdot)$ is actually of class $C^2(\R^N)$. The strategy consists in bounding $u(1,x)$ from above as $|x|\to+\infty$ by a function having the same decay as $u_0$, and then in showing that $u(1+t,\cdot)$ is above $u(1,\cdot)$ in $\R^N$ for all $t>0$ large enough. To do so, we will use some lower and upper bounds for the heat kernel associated with the linearized equation~\eqref{eqv} below, as well as the spreading property~\eqref{defc}. For the sake of clarity, the two cases -- Gaussian decay for $u_0$ or~\eqref{u0}-- will be treated separately.\par
{\it Case 1: Gaussian decay}. Assume here that $u_0$ has Gaussian decay at infinity, that is, there exists $\beta>0$ such that $u_0(x)=O(e^{-\beta|x|^2})$ as $|x|\to+\infty$. Since $u_0\in L^{\infty}(\R^N;[0,1])$, there is then~$C>0$ such that
$$0\le u_0(x)\le C\,e^{-\beta|x|^2}\ \hbox{ for a.e. }x\in\R^N.$$
Remember that the function $f$ is globally Lipschitz continuous in its second variable, uniformly with respect to $x\in\R^N$. Since $f(\cdot,0)=0$ in $\R^N$, let then $L>0$ be such that
\be\label{defL}
f(x,s)\le Ls\ \hbox{ for all }(x,s)\in\R^N\times[0,1].
\ee
The maximum principle yields
$$0\le u(1,x)\le e^L\,v(1,x)\ \hbox{ for all }x\in\R^N,$$
where $v$ denotes the solution of the Cauchy problem
\be\label{eqv}\left\{\baa{rcll}
v_t & = & \div(A(x)\nabla v), & t>0,\ x\in\R^N,\vspace{3pt}\\
v(0,\cdot) & = & u_0.\eaa\right.
\ee
Therefore,
$$0\le u(1,x)\le C\,e^L\,\int_{\R^N}p(1,x;y)\,e^{-\beta|y|^2}\,dy\ \hbox{ for all }x\in\R^N,$$
where $p(t,x;y)$ denotes the heat kernel associated to the linear equation~\eqref{eqv}, that is, for every $y\in\R^N$, $p(\cdot,\cdot;y)$ solves~\eqref{eqv} with the Dirac distribution $\delta_y$ at $y$ as initial condition. It follows from the bounds of $p$ in~\cite{n} (see also~\cite{a,d,fs,fr} for related results) that there is a real number $K\ge1$ such that
\be\label{boundsp}
\frac{e^{-K|x-y|^2/t}}{K\,t^{N/2}}\le p(t,x;y)\le\frac{K\,e^{-|x-y|^2/(Kt)}}{t^{N/2}}\ \hbox{ for all }t>0\hbox{ and }(x,y)\in\R^N\times\R^N.
\ee
In particular,
$$0\le u(1,x)\le K\,C\,e^L\int_{\R^N}e^{-|x-y|^2/K-\beta|y|^2}dy$$
Let $\eta\in(0,1)$ be such that $\eta<\beta\,K\,(1-\eta)$ and denote $\rho=\beta-\eta/(K(1-\eta))>0$. By writing
$$\baa{rcl}
\displaystyle-\frac{|x-y|^2}{K}=-\frac{|x|^2}{K}+\frac{2(x\cdot y)}{K}-\frac{|y|^2}{K} & \le & \displaystyle-\frac{|x|^2}{K}+\frac{(1-\eta)|x|^2}{K}+\frac{|y|^2}{K(1-\eta)}-\frac{|y|^2}{K}\vspace{3pt}\\
& = & \displaystyle-\frac{\eta\,|x|^2}{K}+\frac{\eta\,|y|^2}{K\,(1-\eta)},\eaa$$
it follows that
$$0\le u(1,x)\le K\,C\,e^L\,e^{-\eta|x|^2/K}\int_{\R^N}e^{-\rho|y|^2}dy\ \hbox{ for all }x\in\R^N.$$
To sum up, since the continuous function $u(1,\cdot)$ is less than $1$ in $\R^N$ by~\eqref{0u1}, one infers that there exist some real numbers $\theta\in(0,1)$ and $\omega>0$ such that
\be\label{u1}
u(1,x)\le\min\big(\theta,\omega\,e^{-\eta|x|^2/K}\big)\ \hbox{ for all }x\in\R^N.
\ee\par
Let us now show that $u(1+t,\cdot)$ is above $u(1,\cdot)$ in $\R^N$ for all $t>0$ large enough. Since $f$ is nonnegative in $\R^N\times[0,1]$, one infers from the maximum principle that $u(1+t,x)\ge v(1+t,x)$ for all $t\ge0$ and $x\in\R^N$, where $v$ solves~\eqref{eqv}. Since $u_0$ is nonnegative a.e. in $\R^N$ and non-trivial, there is $R>0$ such that
$$\sigma:=\int_{B(0,R)}u_0(y)\,dy>0$$
and
$$u(1+t,x)\ge\int_{\R^N}p(1+t,x;y)\,u_0(y)\,dy\ge\frac{1}{K\,(1+t)^{N/2}}\int_{B(0,R)}e^{-K|x-y|^2/(1+t)}\,u_0(y)\,dy$$
for all $t\ge0$ and $x\in\R^N$, from~\eqref{boundsp}. By writing
$$-\frac{K\,|x-y|^2}{1+t}\ge-\frac{2\,K\,|x|^2}{1+t}-\frac{2\,K\,|y|^2}{1+t}\ge-\frac{2\,K\,|x|^2}{1+t}-2KR^2$$
for all $t\ge0$, $x\in\R^N$ and $y\in B(0,R)$, one gets that
\be\label{ut1}
u(1+t,x)\ge\frac{e^{-2KR^2}\,e^{-2K|x|^2/(1+t)}}{K\,(1+t)^{N/2}}\int_{B(0,R)}u_0(y)\,dy=\frac{\sigma\,e^{-2KR^2}\,e^{-2K|x|^2/(1+t)}}{K\,(1+t)^{N/2}}
\ee
for all $t\ge0$ and $x\in\R^N$.\par
We finally show that~\eqref{defT} holds for some $T>0$ large enough. Assume not. Then there exist a sequence $(T_n)_{n\in\N}$ of positive real numbers and a sequence $(x_n)_{n\in\N}$ of points in $\R^N$ such that~$T_n\to+\infty$ as $n\to+\infty$ and $u(1+T_n,x_n)<u(1,x_n)$ for all $n\in\N$. Since $u(1,\cdot)\le\theta<1$ in~$\R^N$ and $\min_{|x|\le c t}u(t,x)\to1$ as $t\to+\infty$ with $c>0$ by~\eqref{defc}, it follows that $|x_n|\ge c(1+T_n)$ for $n$ large enough, while $u(1+T_n,x_n)<u(1,x_n)$ and~\eqref{u1}-\eqref{ut1} yield
$$\frac{\sigma\,e^{-2KR^2}\,e^{-2K|x_n|^2/(1+T_n)}}{K\,(1+T_n)^{N/2}}<\omega\,e^{-\eta|x_n|^2/K}\ \hbox{ for all }n\in\N,$$
whence
$$\sigma\,K^{-1}\,\omega^{-1}\,e^{-2KR^2}\,(1+T_n)^{-N/2}<e^{-\eta|x_n|^2/K+2K|x_n|^2/(1+T_n)}\le e^{-\eta|x_n|^2/(2K)}\le e^{-\eta c^2(1+T_n)^2/(2K)}$$
for all $n$ large enough. This clearly leads to a contradiction. As a consequence, there is $T>0$ such that~\eqref{defT} holds.\par
{\it Case 2: assumption~\eqref{u0}}. Since $0\le u_0\le 1$ a.e. in $\R^N$, it follows from~\eqref{u0} that there is~$\delta'>0$ such that $u_0(x)\le\delta'\,e^{-\lambda|x|}$ for a.e. $x\in\R^N$. Therefore, with the same notations as in case~1, one infers that
$$u(1,x)\le\delta'\,e^L\int_{\R^N}p(1,x;y)\,e^{-\lambda|y|}\,dy\le K\,\delta'\,e^L\int_{\R^N}e^{-|x-y|^2/K-\lambda|y|}\,dy\ \hbox{ for all }x\in\R^N.$$
Hence,
$$u(1,x)\le K\,\delta'\,e^L\int_{\R^N}e^{-|y|^2/K-\lambda|x-y|}\,dy\le K\,\delta'\,e^L\,e^{-\lambda|x|}\int_{\R^N}e^{-|y|^2/K+\lambda|y|}\,dy\ \hbox{ for all }x\in\R^N$$
and, since $u(1,\cdot)$ is continuous and less than $1$ in $\R^N$, there are then $\theta'\in(0,1)$ and $\omega'>0$ such that
\be\label{u1bis}
u(1,x)\le\min\big(\theta',\omega'\,e^{-\lambda|x|}\big)\ \hbox{ for all }x\in\R^N.
\ee\par
On the other hand, assumption~\eqref{u0} yields the existence of $R>0$ such that $u_0(x)\ge\gamma\,e^{-\lambda|x|}$ for all $|x|\ge R$. It follows then from~\eqref{boundsp} and the nonnegativity of $f$ and $u_0$ that, for all $t\ge0$ and~$x\in\R^N$,
\be\label{uT1bis}\baa{rcl}
\displaystyle u(1+t,x)\ge\int_{\R^N}p(1+t,x;y)\,u_0(y)\,dy & \!\!\ge\!\! & \displaystyle\frac{\gamma}{K\,(1+t)^{N/2}}\int_{\R^N\backslash B(0,R)}e^{-K|x-y|^2/(1+t)-\lambda|y|}\,dy\vspace{3pt}\\
& \!\!=\!\! & \displaystyle\frac{\gamma}{K}\int_{\{z\in\R^N;\,|x-\sqrt{1+t}\,z|\ge R\}}e^{-K|z|^2-\lambda|x-\sqrt{1+t}\,z|}\,dz.\eaa
\ee\par
Assume now by contradiction that property~\eqref{defT} does not hold for any $T>0$. Then there exist a sequence $(T_n)_{n\in\N}$ of positive real numbers and a sequence $(x_n)_{n\in\N}$ of points in $\R^N$ such that $T_n\to+\infty$ as $n\to+\infty$ and $u(1+T_n,x_n)<u(1,x_n)$ for all $n\in\N$. Since $u(1,\cdot)\le\theta'<1$ in~$\R^N$ and $\min_{|x|\le c t}u(t,x)\to1$ as $t\to+\infty$ with $c>0$ by~\eqref{defc}, it follows that $|x_n|\ge c(1+T_n)$ for $n$ large enough, while $u(1+T_n,x_n)<u(1,x_n)$ and~\eqref{u1bis}-\eqref{uT1bis} yield
$$\omega'\,e^{-\lambda|x_n|}>\frac{\gamma}{K}\int_{\{z\in\R^N;\,|x_n-\sqrt{1+T_n}\,z|\ge R\}}e^{-K|z|^2-\lambda|x_n-\sqrt{1+T_n}\,z|}\,dz\ \hbox{ for all }n\in\N.$$
Since $\liminf_{n\to+\infty}|x_n|/T_n\ge c>0$, one has $B(x_n/|x_n|,1/2)\subset\{z\in\R^N;\,|x_n-\sqrt{1+T_n}\,z|\ge R\}$ for~$n$ large enough, whence
$$\omega'\,e^{-\lambda|x_n|}>\frac{\gamma}{K}\int_{B(x_n/|x_n|,1/2)}\!\!\!e^{-K|z|^2-\lambda|x_n-\sqrt{1+T_n}\,z|}\,dz\ge\frac{\gamma\,e^{-9K/4}}{K}\int_{B(0,1/2)}\!\!\!e^{-\lambda|x_n-\sqrt{1+T_n}\,(x_n/|x_n|+y)|}\,dy$$
for $n$ large enough. For $n$ large enough so that $\sqrt{1+T_n}\le|x_n|$, it follows that, for all~$y\in B(0,1/2)$,
$$\Big|x_n-\sqrt{1+T_n}\,\Big(\frac{x_n}{|x_n|}+y\Big)\Big|\le|x_n|\Big(1-\frac{\sqrt{1+T_n}}{|x_n|}\Big)+\frac{\sqrt{1+T_n}}{2}=|x_n|-\frac{\sqrt{1+T_n}}{2},$$
whence
$$\omega'\,e^{-\lambda|x_n|}>\frac{\gamma\,e^{-9K/4}\,e^{-\lambda|x_n|+\lambda\sqrt{1+T_n}/2}}{K}\int_{B(0,1/2)}dy$$
for $n$ large enough. This leads to a contradiction since $T_n\to+\infty$ as $n\to+\infty$.\par
As a conclusion,~\eqref{defT} holds when~\eqref{u0} is fulfilled and the proof of Lemma~\ref{lem1} is thereby complete.~\hfill$\Box$\break

From Lemma~\ref{lem1} and the maximum principle, the following corollary immediately holds.

\begin{cor}\label{cor1}
For every $t\ge1$, $T'\ge T$ and $x\in\R^N$, one has $u(t+T',x)\ge u(t,x)$.
\end{cor}

\begin{rem}{\rm Notice from the proof of Lemma~\ref{lem1} that time $1$ could be replaced by any positive time $t_0$ in the statement: namely, for any $t_0>0$, there exists $T_0>0$ such that~$u(t_0+t,x)\ge u(t_0,x)$ for all $t\ge T_0$ and $x\in\R^N$. However, this property does not hold in general with $t_0=0$. Indeed, if~$0\le u_0\le 1$ is continuous and $\max_{\R^N}u_0=1$, then $u_0$ can never be bounded from above in $\R^N$ by~$u(t,\cdot)$ for any $t>0$, since $u(t,x)<1$ for all $t>0$ and~$x\in\R^N$ by the strong parabolic maximum principle.}
\end{rem}

\begin{rem}{\rm The assumptions~\eqref{u0gaussian} or~\eqref{u0} were crucially used in the proof of Lemma~\ref{lem1}, in order to trap $u(1,x)$ between two comparable functions as $|x|\to+\infty$, the lower one giving rise to a solution which, after some time, is above the upper one at time $1$. The conclusion of Lemma~\ref{lem1} may not hold for more general initial conditions $u_0$, for instance if $\gamma\,e^{-\lambda_1|x|}\le u_0(x)\le\delta\,e^{-\lambda_2|x|}$ for $|x|$ large enough, with $\gamma,\,\delta>0$, $0<\lambda_2<\lambda_1$, $\liminf_{|x|\to+\infty}u_0(x)\,e^{\lambda_1|x|}<+\infty$ and $\limsup_{|x|\to+\infty}u_0(x)\,e^{\lambda_2|x|}>0$. For such initial conditions, more complex dynamics may occur in general, even for homogeneous one-dimensional equations, see e.g.~\cite{hn,y}.}
\end{rem}


\SE{Improved monotonicity when $u(t,x)$ is away from $0$ and~$1$}\label{sec3}

In this section, we improve the $T$-monotonicity result stated in Corollary~\ref{cor1}, for the points~$(t,x)$ such that $0<a\le u(t,x)\le b<1$, where $0<a\le b<1$ are given. To do so, let us first define
\be\label{deftau*}
\tau_*=\inf\big\{\tau>0;\ \exists\,t_0\ge 0,\ \forall\,\tau'\ge\tau,\ \forall\,t\ge t_0,\ \forall\,x\in\R^N,\ u(t+\tau',x)\ge u(t,x)\big\}.
\ee
It follows from Corollary~\ref{cor1} that $0\le \tau_*\le T<+\infty$. Our goal is to show that $\tau_*=0$ (this goal will be achieved at the beginning of Section~\ref{sec7}). 

\begin{lem}\label{lem2}
Let $a$ and $b$ be any two real numbers such that $0<a\le b<1$ and let $\tau$ be any real number such that $\tau\ge\tau_*$ and $\tau>0$. Then,
$$\liminf_{t\to+\infty,\ a\le u(t,x)\le b}\frac{u(t+\tau,x)}{u(t,x)}>1,$$
that is, there exist $t_0>0$ and $\delta>0$ such that, for all $(t,x)\in[t_0,+\infty)\times\R^N$ with $a\le u(t,x)\le b$, there holds $u(t+\tau,x)\ge(1+\delta)\,u(t,x)$.
\end{lem}

\noindent{\bf{Proof.}} The proof shall use the definition of $\tau_*$ and the positivity of $\tau$ together with the spreading pro\-perties of solutions of equations obtained as finite or infinite spatial shifts of~\eqref{eq}. We argue by contradiction. So, assume that the conclusion of Lemma~\ref{lem2} does not hold. Then there are two sequences $(t_n)_{n\in\N}$ and $(\delta_n)_{n\in\N}$ of positive real numbers and a sequence $(x_n)_{n\in\N}$ of points in $\R^N$ such that $\delta_n\to0$ as $n\to+\infty$, $t_n\to+\infty$ as $n\to+\infty$ and
\be\label{un}
a\le u(t_n,x_n)\le b\ \hbox{ and }\ u(t_n+\tau,x_n)<(1+\delta_n)\,u(t_n,x_n)\ \hbox{ for all }n\in\N.
\ee\par
Shift the origin at the points $(t_n,x_n)$ and define
$$u_n(t,x)=u(t+t_n,x+x_n).$$
The functions $u_n$ are classical solutions of
\be\label{equn}
(u_n)_t=\hbox{div}(A(x+x_n)\nabla u_n)+f(x+x_n,u_n),\ \ t>-t_n,\ x\in\R^N
\ee
with $0<u_n(t,x)<1$ for all $(t,x)\in(-t_n,+\infty)\times\R^N$. From Arzela-Ascoli theorem, up to extraction of a subsequence, the functions $\R^N\times[0,1]\ni(x,s)\mapsto f(x+x_n,s)$ converge locally uniformly in $\R^N\times[0,1]$ to a continuous function $f_{\infty}:\R^N\times[0,1]$ which actually shares with $f$ the following properties: $f_{\infty}(\cdot,0)=f_{\infty}(\cdot,1)=0$, $f_{\infty}(x,1-u)/u$ is nonincreasing in $u\in(0,1]$, and $f_{\infty}$ satisfies~\eqref{hypf}, whence $\inf_{x\in\R^N}f_{\infty}(x,s)>0$ for every $s\in(0,1)$). Furthermore, up to extraction of another subsequence, the matrix fields $x\mapsto A(x+x_n)$ converge in $C^1_{loc}(\R^N)$ to a uniformly definite positive symmetric matrix field $A_{\infty}$.\footnote{As a matter of fact, since $u(t_n,x_n)\le b<1$ and $t_n\to+\infty$, then $|x_n|\to+\infty$ by~\eqref{defc}, whence $A_{\infty}$ is a constant matrix due to~\eqref{oscA}. However, the fact that $A_{\infty}$ is constant is not used in the proof of the present lemma.} Lastly, from standard parabolic estimates, the functions $u_n$ converge locally uniformly in $C^{1,2}_{t,x}(\R\times\R^N)$, up to extraction of another subsequence, to a classical solution~$u_{\infty}$ of
\be\label{uinfty}
(u_{\infty})_t=\hbox{div}(A_{\infty}\nabla u_{\infty})+f_{\infty}(x,u_{\infty}),\ \ t\in\R,\ x\in\R^N,
\ee
such that $0\le u_{\infty}(t,x)\le 1$ for all $(t,x)\in\R\times\R^N$.\par
Now, for any $\epsilon>0$, it follows from $\tau\ge\tau_*$ and from the definition of~$\tau_*$ in~\eqref{deftau*} that there is~$T_0>0$ such that
$$u(t+\tau+\epsilon,x)\ge u(t,x)\ \hbox{ for all }(t,x)\in[T_0,+\infty)\times\R^N$$
(actually, if $\tau>\tau_*$, then one can also take $\epsilon=0$). In particular, since $t_n\to+\infty$ as $n\to+\infty$, one infers that $u_{\infty}(t+\tau+\epsilon,x)\ge u_{\infty}(t,x)$ for all $(t,x)\in\R\times\R^N$. Since $\epsilon>0$ can be arbitrary, one gets that
$$u_{\infty}(t+\tau,x)\ge u_{\infty}(t,x)\ \hbox{ for all }(t,x)\in\R\times\R^N.$$
On the other hand, the inequalities~\eqref{un} and $\lim_{n\to+\infty}\delta_n=0$ imply that $a\le u_{\infty}(0,0)\le b$ and $u_{\infty}(\tau,0)\le u_{\infty}(0,0)$, whence $u_{\infty}(\tau,0)=u_{\infty}(0,0)$. As a consequence, the bounded functions~$u_{\infty}(\cdot+\tau,\cdot)$ and $u_{\infty}(\cdot,\cdot)$ are ordered in $\R\times\R^N$ and are equal at $(0,0)$. It follows from the strong maximum principle that
$$u_{\infty}(t+\tau,x)=u_{\infty}(t,x)$$
for all $(t,x)\in(-\infty,0]\times\R^N$, and then for all $(t,x)\in\R\times\R^N$ from the uniqueness of the Cauchy problem associated with~\eqref{uinfty}. Furthermore, $0<a\le u_{\infty}(0,0)\le b<1$ and~$0\le u_{\infty}\le1$ in~$\R\times\R^N$, whence $0<u_{\infty}<1$ in $\R\times\R^N$ from the strong maximum principle. Lastly, $u_{\infty}(t,x)\to1$ as~$t\to+\infty$ locally uniformly in $x\in\R^N$, as recalled in Section~\ref{intro} for $u$ and $f$, from the properties shared by $f_{\infty}$ with $f$. Thus, the limit $\N\ni m\to+\infty$ in $u_{\infty}(m\tau,0)=u_{\infty}(0,0)\le b<1$ leads to a contradiction, since $\tau>0$ by assumption. The proof of Lemma~\ref{lem2} is thereby complete.~\hfill$\Box$\break

From Lemma~\ref{lem2} and the uniform continuity of $u$ in, say, $[1,+\infty)\times\R^N$, the inequalities stated in Lemma~\ref{lem2} hold uniformly for some time-shifts in a neighborhood of $\tau_*$ if $\tau_*$ is positive, as the following corollary shows.

\begin{cor}\label{cor2}
Let $a$ and $b$ be any two real numbers such that $0<a\le b<1$. If one assumes that $\tau^*>0$, then there exist~$t_0>0$, $\delta>0$ and $0<\underline{\tau}<\tau_*<\overline{\tau}$ such that, for all $\tau\in[\underline{\tau},\overline{\tau}]$ and $(t,x)\in[t_0,+\infty)\times\R^N$ with~$a\le u(t,x)\le b$, then $u(t+\tau,x)\ge(1+\delta)\,u(t,x)$.
\end{cor}

\noindent{\bf{Proof.}} From Lemma~\ref{lem2} applied with $\tau=\tau_*$, there are $t_0>0$ and $\delta>0$ such that~$u(t+\tau_*,x)\ge(1+2\delta)\,u(t,x)$ for all $(t,x)\in[t_0,+\infty)\times\R^N$ with $a\le u(t,x)\le b$. Choose~$\epsilon\in(0,1)$ so that $(1-\epsilon)(1+2\delta)\ge1+\delta$. Since $u$ is uniformly continuous in~$[t_0,+\infty)\times\R^N$ from standard parabolic estimates, there exist some real numbers~$\underline{\tau}$ and~$\overline{\tau}$ such that $0<\underline{\tau}<\tau_*<\overline{\tau}$ and
$$|u(t+\tau,x)-u(t+\tau_*,x)|\le\epsilon\,(1+2\delta)\,a\ \hbox{ for all }\tau\in[\underline{\tau},\overline{\tau}]\hbox{ and for all }(t,x)\in[t_0,+\infty)\times\R^N.$$
Fix now any $\tau\in[\underline{\tau},\overline{\tau}]$ and any $(t,x)\in[t_0,+\infty)\times\R^N$ with $a\le u(t,x)\le b$. One has~$u(t+\tau_*,x)\ge(1+2\delta)\,u(t,x)\ge(1+2\delta)\,a$, whence
$$u(t+\tau,x)\ge u(t+\tau_*,x)-\epsilon\,(1+2\delta)\,a\ge(1-\epsilon)\,u(t+\tau_*,x)\ge(1-\epsilon)\,(1+2\delta)\,u(t,x)\ge(1+\delta)\,u(t,x).$$
This is the desired result and the proof is thereby complete.\hfill$\Box$


\SE{Improved monotonicity when $u(t,x)$ is away from $0$}\label{sec4}

In this section, by using especially the fact that $f(x,1-u)/u$ is nonincreasing with respect to $u\in(0,1]$ for every $x\in\R^N$, we improve the $\tau$-monotonicity of $u$ (with $\tau>\tau_*$) in the region where $u(t,x)$ is close to $1$ (we recall that~$0<u(t,x)<1$ for all~$(t,x)\in(0,+\infty)\times\R^N$). Namely, we will prove the following lemma.

\begin{lem}\label{lem3}
Let $a$ and $\tau$ be any real numbers such that $0<a<1$ and $\tau>\tau_*$. Then,
$$\limsup_{t\to+\infty,\ u(t,x)\ge a}\frac{1-u(t+\tau,x)}{1-u(t,x)}<1,$$
that is, there exist $t_0>0$ and $\delta>0$ such that, for all $(t,x)\in[t_0,+\infty)\times\R^N$ with $u(t,x)\ge a$, there holds $1-u(t+\tau,x)\le(1-\delta)\,(1-u(t,x))$.
\end{lem}

\noindent{\bf{Proof.}} First of all, since $\tau>\tau_*$, it follows from the definition of $\tau_*$ that there is $T_0>0$ such that
\be\label{utau}
u(t+\tau,x)\ge u(t,x)\hbox{ for all }(t,x)\in[T_0,+\infty)\times\R^N.
\ee
Notice that the strong maximum principle then yields $u(t+\tau,x)>u(t,x)$ in $(T_0,+\infty)\times\R^N$ (otherwise, one would have $u(t+\tau,x)=u(t,x)$ in $[T_0,T_1]\times\R^N$ with some $T_1>T_0$, whence~$u(t+\tau,x)=u(t,x)$ in $[T_0,+\infty)\times\R^N$ and $u(T_0+m\tau,0)=u(T_0,0)<1$ for all $m\in\N$, whereas~$u(t,0)\to1$ as $t\to+\infty$. Even if it means increasing $T_0$, one can then assume without loss of generality that
$$u(t+\tau,x)>u(t,x)\hbox{ for all }(t,x)\in[T_0,+\infty)\times\R^N,\ \ u(T_0,0)\ge a\ \hbox{ and }\ T_0>\tau.$$\par
Define now, for every $k\in\N=\{0,1,2,\cdots\}$,
$$E_k=\big\{x\in\R^N;\,\exists\,t\in[T_0+k\tau,T_0+(k+1)\tau],\,u(t,x)\ge a\big\}.$$
The set $E_0$ is not empty since $u(T_0,0)\ge a$. As a consequence,
$$u(T_0+k\tau,0)\ge u(T_0+(k-1)\tau,0)\ge\cdots\ge u(T_0,0)\ge a,$$
whence $0\in E_k$ for every $k\in\N$. Thanks to~\eqref{utau}, the same argument implies that $E_k\subset E_{k+1}$ for every $k\in\N$. Furthermore, each set $E_k$ is closed by continuity of $u$ in $[T_0,+\infty)\times\R^N$. Lastly, as done for the proof of~\eqref{u1} and~\eqref{u1bis} in Lemma~\ref{lem1}, one easily infers that $u(t,x)\to0$ as~$|x|\to+\infty$ locally uniformly in $t>0$, whence each set $E_k$ is bounded. Therefore, the sets~$E_k$ are a non-decreasing sequence of non-empty compact subsets of $\R^N$.\par
We are going to apply the maximum principle to the functions $1-u(t+\tau,x)$ and $1-u(t,x)$ in the sets $[T_0+k\tau,T_0+(k+1)\tau]\times E_k$ by induction with respect to $k$, in order to improve quantitatively the inequality $1-u(t+\tau,x)\le 1-u(t,x)$ in $[T_0+k\tau,T_0+(k+1)\tau]\times E_k$.\par
To do so, we first claim that the function $u$ is bounded from below by a positive constant uniformly in the sets $[T_0+k\tau,T_0+(k+1)\tau]\times E_k$, that is, there is $\underline{a}\in(0,a]$ such that
\be\label{claim1}
\forall\,k\in\N,\ \forall\,(t,x)\in[T_0+k\tau,T_0+(k+1)\tau]\times E_k,\ u(t,x)\ge\underline{a}>0.
\ee
Indeed, otherwise, there exist a sequence $(k_n)_{n\in\N}$ of integers and, for each $n\in\N$, a time~$t_n\in[T_0+k_n\tau,T_0+(k_n+1)\tau]$ and a point $x_n\in E_{k_n}$, with $u(t_n,x_n)\to0$ as $n\to+\infty$. For each $n\in\N$, since $x_n\in E_{k_n}$, there is a time $t'_n\in[T_0+k_n\tau,T_0+(k_n+1)\tau]$ such that $u(t'_n,x_n)\ge a$. Consider the functions
$$(t,x)\mapsto u_n(t,x)=u(t+t_n,x+x_n),$$
which are defined in $(-t_n,+\infty)\times\R^N\supset(-T_0,+\infty)\times\R^N$ and solve~\eqref{equn}, together with $0\le u_n\le 1$. From Arzela-Ascoli theorem and standard parabolic estimates, up to extraction of a subsequence, these functions $u_n$ converge locally uniformly in $C^{1,2}_{t,x}((-T_0,+\infty)\times\R^N)$ to a solution $0\le u_{\infty}\le 1$ of an equation of the type~\eqref{uinfty} in $(-T_0,+\infty)\times\R^N$ (notice that the sequences $(t_n)_{n\in\N}$ and $(x_n)_{n\in\N}$ may not be unbounded and the limiting equation satisfied by~$u_{\infty}$ may just be a finite spatial shift of~\eqref{eq}). Anyway, $u_n(0,0)=u(t_n,x_n)\to0$ as $n\to+\infty$, whence $u_{\infty}(0,0)=0$. Therefore, $u_{\infty}=0$ in $(-T_0,0]\times\R^N$ from the strong maximum principle, and $u_{\infty}=0$ in $(-T_0,+\infty)\times\R^N$ from the uniqueness of the Cauchy problem associated with~\eqref{uinfty}. On the other hand, $|t'_n-t_n|\le\tau<T_0$ for every $n\in\N$. Up to extraction of another subsequence, one can assume that $t'_n-t_n\to t'_{\infty}>-T_0$ as $n\to+\infty$. Since $u_n(t'_n-t_n,0)=u(t'_n,x_n)\ge a$, one gets $u_{\infty}(t'_{\infty},0)\ge a>0$, which leads to a contradiction. As a consequence, the claim~\eqref{claim1} is proved.\par
The second claim is concerned with an upper bound of the values of $u$ on the boundaries~$\partial E_k$ of the sets $E_k$, on the time intervals $[T_0+k\tau,T_0+(k+1)\tau]$. Namely, we claim that there is a real number $b\in(0,1)$ such that
\be\label{claim2}
\forall\,k\in\N,\ \forall\,(t,x)\in[T_0+k\tau,T_0+(k+1)\tau]\times\partial E_k,\ u(t,x)\le b<1.
\ee
Assume not. Then, there exist a sequence $(k_n)_{n\in\N}$ of integers and, for each $n\in\N$, a time~$t_n\in[T_0+k_n\tau,T_0+(k_n+1)\tau]$ and a point~$x_n\in\partial E_{k_n}$, with $u(t_n,x_n)\to1$ as $n\to+\infty$. For each~$n\in\N$, since $x_n\in\partial E_{k_n}\subset E_{k_n}$ and since $u(\cdot,x_n)$ is continuous on $[T_0+k_n\tau,T_0+(k_n+1)\tau]$, the definition of $E_{k_n}$ yields
$$\max_{[T_0+k_n\tau,T_0+(k_n+1)\tau]}u(\cdot,x_n)\ge a.$$
Furthermore, if $\min_{[T_0+k_n\tau,T_0+(k_n+1)\tau]}u(\cdot,x_n)>a$, then by uniform continuity of $u$ in $[T_0,+\infty)\times\R^N$ one would have $\min_{[T_0+k_n\tau,T_0+(k_n+1)\tau]}u(\cdot,x)>a$ for all $x$ in a neighborhood of~$x_n$ and~$x_n$ would then be an interior point of $E_{k_n}$. Therefore, $\min_{[T_0+k_n\tau,T_0+(k_n+1)\tau]}u(\cdot,x_n)\le a$ and there is a time~$t'_n\in[T_0+k_n\tau,T_0+(k_n+1)\tau]$ such that
$$u(t'_n,x_n)=a.$$
Now, as in the previous paragraph, the functions $(t,x)\mapsto u_n(t,x)=u(t+t_n,x+x_n)$ converge, up to extraction of a subsequence, locally uniformly in $C^{1,2}_{t,x}((-T_0,+\infty)\times\R^N)$ to a solution $0\le u_{\infty}\le 1$ of an equation of the type~\eqref{uinfty} in $(-T_0,+\infty)\times\R^N$. One has $u_{\infty}(0,0)=1$, whence $u_{\infty}=1$ in $(-T_0,0]\times\R^N$ and then in $(-T_0,+\infty)\times\R^N$. On the other hand, up to extraction of another subsequence, there holds $\lim_{n\to+\infty}(t'_n-t_n)=t'_{\infty}\in[-\tau,\tau]\subset(-T_0,+\infty)$ and $u_{\infty}(t'_{\infty},0)=a<1$. One has reached a contradiction, and the claim~\eqref{claim2} follows.\par
Similarly, we claim that there is a real number $\overline{b}\in(0,1)$ such that
\be\label{claim3}
\forall\,k\in\N,\ \forall\,x\in E_{k+1}\backslash E_k,\ u(T_0+(k+1)\tau,x)\le\overline{b}<1.
\ee
Otherwise, there exist a sequence $(k_n)_{n\in\N}$ of integers and, for each $n\in\N$, a point $x_n\in E_{k_n+1}\backslash E_{k_n}$, with $u(T_0+(k_n+1)\tau,x_n)\to1$ as $n\to+\infty$. For each $n\in\N$, since $x_n\in E_{k_n+1}\backslash E_{k_n}$, there holds
$$\max_{[T_0+(k_n+1)\tau,T_0+(k_n+2)\tau]}u(\cdot,x_n)\ge a\ \hbox{ and }\ \max_{[T_0+k_n\tau,T_0+(k_n+1)\tau]}u(\cdot,x_n)<a,$$
whence there is a time $t_n\in[T_0+(k_n+1)\tau,T_0+(k_n+2)\tau]$ such that $u(t_n,x_n)=a$. Up to extraction of a subsequence, the functions
$$(t,x)\mapsto u_n(t,x)=u(t+T_0+(k_n+1)\tau,x+x_n)$$
converge locally uniformly in $C^{1,2}_{t,x}((-T_0-\tau,+\infty)\times\R^N)$ to a solution $0\le u_{\infty}\le 1$ of an equation of the type~\eqref{uinfty} in $(-T_0-\tau,+\infty)\times\R^N$. One has $u_{\infty}(0,0)=1$, whence $u_{\infty}=1$ in~$(-T_0-\tau,0]\times\R^N$ and then in $(-T_0-\tau,+\infty)\times\R^N$. On the other hand, up to extraction of another subsequence, there holds $\lim_{n\to+\infty}(t_n-(T_0+(k_n+1)\tau))=t_{\infty}\in[0,\tau]\subset(-T_0-\tau,+\infty)$ and $u_{\infty}(t_{\infty},0)=a<1$. One has reached a contradiction, and the claim~\eqref{claim3} is proved.\par
Putting together~\eqref{claim1},~\eqref{claim2} and~\eqref{claim3}, one gets that
$$\left\{\baa{ll}
\forall\,k\in\N,\ \forall\,(t,x)\in[T_0+k\tau,T_0+(k+1)\tau]\times\partial E_k, & 0<\underline{a}\le u(t,x)\le b<1,\vspace{3pt}\\
\forall\,k\in\N,\ \forall\,x\in E_{k+1}\backslash E_k, & 0<\underline{a}\le u(T_0+(k+1)\tau,x)\le\overline{b}<1.\eaa\right.$$
It follows then from Lemma~\ref{lem2} applied once with $(\underline{a},b,\tau)$ and another time with $(\underline{a},\overline{b},\tau)$ (notice that $\tau\ge\tau_*$ and $\tau>0$ since here $\tau>\tau_*$) that there are $k_0\in\N$ and $\delta_0\in(0,+\infty)$ such that, for all $k\ge k_0$,
$$\left\{\baa{ll}
\forall\,(t,x)\in[T_0+k\tau,T_0+(k+1)\tau]\times\partial E_k, & u(t+\tau,x)\ge(1+\delta_0)\,u(t,x),\vspace{3pt}\\
\forall\,x\in E_{k+1}\backslash E_k, & u(T_0+(k+2)\tau,x)\ge(1+\delta_0)\,u(T_0+(k+1)\tau,x).\eaa\right.$$
Define
$$\delta=\delta_0\underline{a}>0.$$
One infers that
\be\label{Ek}\baa{l}
\forall\,k\ge k_0,\ \forall\,(t,x)\in[T_0+k\tau,T_0+(k+1)\tau]\times\partial E_k,\vspace{3pt}\\
\quad\baa{rcl}
1-u(t+\tau,x)\le1-(1+\delta_0)\,u(t,x)\le1-u(t,x)-\delta_0\underline{a} & \!\!=\!\! & 1-u(t,x)-\delta\vspace{3pt}\\
& \!\!\le\!\! & (1-\delta)\,(1-u(t,x)),\eaa\eaa
\ee
and, by arguing similarly with $x\in E_{k+1}\backslash E_k$, that
\be\label{Ek1k}
\forall\,k\ge k_0,\ \forall\,x\in E_{k+1}\backslash E_k,\ \ 1-u(T_0+(k+2)\tau,x)\le(1-\delta)\,(1-u(T_0+(k+1)\tau,x)).
\ee
On the other hand, since
$$1>u(t+\tau,x)>u(t,x)>0$$
for all $(t,x)\in[T_0+k_0\tau,T_0+(k_0+1)\tau]\times E_{k_0}$ and since both functions $u(\cdot+\tau,\cdot)$ and $u$ are continuous on this compact set $[T_0+k_0\tau,T_0+(k_0+1)\tau]\times E_{k_0}$, if follows that, even if it means decreasing $\delta>0$,
\be\label{Ek0}
\forall\,(t,x)\in[T_0+k_0\tau,T_0+(k_0+1)\tau]\times E_{k_0},\ \ 1-u(t+\tau,x)\le(1-\delta)\,(1-u(t,x)).
\ee\par
Finally, we claim by induction on $k$ that
\be\label{claim4}
\forall\,k\ge k_0,\ \forall\,(t,x)\in[T_0+k\tau,T_0+(k+1)\tau]\times E_k,\ \ 1-u(t+\tau,x)\le(1-\delta)\,(1-u(t,x)).
\ee
First of all, the property is true at $k=k_0$, by~\eqref{Ek0}. Assume now that the property is satisfied for some $k\in\N$ with $k\ge k_0$. In particular, by choosing $t=T_0+(k+1)\tau$, there holds
$$\forall\,x\in E_k,\ \ 1-u(t_0+(k+2)\tau,x)\le(1-\delta)\,(1-u(t_0+(k+1)\tau,x)).$$
This last inequality also holds for all $x\in E_{k+1}\backslash E_k$, by~\eqref{Ek1k}. Therefore,
\be\label{Ek1}
\forall\,x\in E_{k+1},\ \ 1-u(t_0+(k+2)\tau,x)\le(1-\delta)\,(1-u(t_0+(k+1)\tau,x)).
\ee
Furthermore, property~\eqref{Ek} yields
\be\label{Ek1bis}
\forall\,(t,x)\in[t_0+(k+1)\tau,t_0+(k+2)\tau]\times\partial E_{k+1},\ \ 1-u(t+\tau,x)\le(1-\delta)\,(1-u(t,x)).
\ee
Consider the functions
$$v(t,x)=1-u(t+\tau,x)\ \hbox{ and }\ \overline{v}(t,x)=(1-\delta)\,(1-u(t,x))$$
in the compact set
$$Q_k=[t_0+(k+1)\tau,t_0+(k+2)\tau]\times E_{k+1}.$$
The inequalities~\eqref{Ek1} and~\eqref{Ek1bis} mean that
$$v(t,x)\le\overline{v}(t,x)\ \hbox{ for all }(t,x)\,\in\,\{t_0+(k+1)\tau\}\!\times\!E_{k+1}\,\cup\,[t_0+(k+1)\tau,t_0+(k+2)\tau]\!\times\!\partial E_{k+1},$$
namely $v\le\overline{v}$ on the parabolic boundary of $Q_k$. Let us now check that $\overline{v}$ is a supersolution of the equation satisfied by $v$. On the one hand, the function $v$ satisfies $0\le v\le 1$ and obeys
$$v_t=\hbox{div}(A(x)\nabla v)+g(x,v)\ \hbox{ in }Q_k,$$
where $g$ is defined by $g(x,s)=-f(x,1-s)$ for all $(x,s)\in\R^N\times[0,1]$. On the other hand, the function $\overline{v}$ satisfies $0\le\overline{v}\le 1$ in $Q_k$ and
$$\baa{rcl}
\overline{v}_t-\hbox{div}(A(x)\nabla\overline{v})-g(x,\overline{v}) & = & -(1-\delta)\,u_t+(1-\delta)\,\hbox{div}(A(x)\nabla u)-g(x,\overline{v})\vspace{3pt}\\
& = & -(1-\delta)\,f(x,u)-g(x,\overline{v})\vspace{3pt}\\
& = & (1-\delta)\,g(x,1-u)-g(x,(1-\delta)\,(1-u)).\eaa$$
But the function $g(x,s)/s$ is nondecreasing with respect to $s\in(0,1]$, since by assumption the function $f(x,1-s)/s$ is nonincreasing with respect to $s\in(0,1]$. Hence,
$$g(x,(1-\delta)\,(1-u(t,x)))\le(1-\delta)\,g(x,1-u(t,x))\ \hbox{ in }Q_k$$
and
$$\overline{v}_t-\hbox{div}(A(x)\nabla\overline{v})-g(x,\overline{v})\ge0\ \hbox{ in }Q_k.$$
The parabolic maximum principle then implies that $v\le\overline{v}$ in $Q_k$. This means that property~\eqref{claim4} is satisfied with $k+1$ and finally that it holds by induction for all $k\ge k_0$.\par
As a conclusion, set $t_0=T_0+k_0\tau$ and consider any $(t,x)\in[t_0,+\infty)\times\R^N$ such that $u(t,x)\ge a$. Let $k\in\N$, $k\ge k_0$ be such that $T_0+k\tau\le t\le T_0+(k+1)\tau$. Thus, $x\in E_k$ and property~\eqref{claim4} yields
$$1-u(t+\tau,x)\le(1-\delta)\,(1-u(t,x)).$$
The proof of Lemma~\ref{lem3} is thereby complete.\hfill$\Box$


\SE{Monotonicity in time when $u(t,x)$ is close to $1$}\label{sec5}

In this section, based on Lemma~\ref{lem3}, we will show that $u$ is actually increasing in time at large time when it is close to $1$.

\begin{lem}\label{lem4}
There exist $b\in(0,1)$ and $\tilde{T}>0$ such that, for all $(t,x)\in[\tilde{T},+\infty)\times\R^N$ with~$u(t,x)\ge b$, there holds $u_t(t,x)>0$.
\end{lem}

\noindent{\bf{Proof.}} As in the proof of Lemma~\ref{lem3}, denote $v=1-u$. The function $v$ satisfies $0<v<1$ in~$(0,+\infty)\times\R^N$ and
$$v_t=\hbox{div}(A(x)\nabla v)+g(x,v),\ \ t>0,\ x\in\R^N$$
with $g(x,s)=-f(x,1-s)$. Furthermore, by choosing, say, $\tau=\tau_*+1$, it follows from definition~\eqref{deftau*} that there is $t_0>1$ such that, for all $(t,x)\in[t_0,+\infty)\times\R^N$, $1>u(t+\tau,x)\ge u(t,x)$, that is,
\be\label{vtau}
0<v(t+\tau,x)\le v(t,x)\ \hbox{ in }[t_0,+\infty)\times\R^N.
\ee
From standard parabolic estimates and Harnack inequality, there are some positive constants $C_1$ and $C_2$ such that
$$\forall\,(t,x)\in[t_0,+\infty)\times\R^N,\ |v_t(t,x)|+|\nabla v(t,x)|\le C_1\max_{[t-1,t]\times\overline{B(x,1)}}v\le C_2\,v(t+\tau,x).$$
Together with~\eqref{vtau}, it follows that the fields $v_t/v$ and $\nabla v/v$ are bounded in $[t_0,+\infty)\times\R^N$. Define now
\be\label{defM}
M=\limsup_{t\to+\infty,\,v(t,x)\to0}\frac{v_t(t,x)}{v(t,x)}.
\ee
From the previous observations and the fact that $v(t,x)=1-u(t,x)\to0$ as $t\to+\infty$ locally uniformly in $x\in\R^N$, one infers that $M$ is a real number. To complete the proof of Lemma~\ref{lem4}, it will actually be sufficient to show that $M<0$.\par
To do so, owing to the definition of $M$, pick a sequence of points $(t_n,x_n)_{n\in\N}$ in $[t_0,+\infty)\times\R^N$ such that
$$t_n\to+\infty,\ \ v(t_n,x_n)\to0\ \hbox{ and }\ \frac{v_t(t_n,x_n)}{v(t_n,x_n)}\to M\ \hbox{ as }n\to+\infty.$$
Define
$$v_n(t,x)=\frac{v(t+t_n,x+x_n)}{v(t_n,x_n)}>0\ \hbox{ in }(-t_n,+\infty)\times\R^N.$$
Since the fields $v_t/v$ and $\nabla v/v$ are bounded in $[t_0,+\infty)\times\R^N$, one infers that the functions $v_n$ are bounded locally in $\R\times\R^N$, in the sense that, for any compact subset $K$ of $\R\times\R^N$, there is $n_K\in\N$ such that $v_n$ is well defined in $K$ for every $n\ge n_K$ and $\sup_{n\ge n_K}\|v_n\|_{L^{\infty}(K)}<+\infty$. Furthermore, the functions $v_n$ obey
\be\label{eqvn}
(v_n)_t(t,x)=\hbox{div}(A(x+x_n)\nabla v_n(t,x))+\frac{g(x+x_n,v(t_n,x_n)\,v_n(t,x))}{v(t_n,x_n)},\ \ t>-t_n,\ x\in\R^N.
\ee
Remember now that $f(\cdot,1)=0$ in $\R^N$, that the function $(x,s)\mapsto f(x,s)$ is Lipschitz continuous with respect to $s$ uniformly in $x\in\R^N$, of class $C^1$ with respect to $s$ in $\R^N\times[s_1,1]$ for some~$s_1\in(0,1)$, and that $f_s$ is uniformly continuous in $\R^N\times[s_1,1]$ and of class $C^{0,\alpha}$ with respect to $x$ uniformly in $s\in[s_1,1]$. Therefore, the function $g$ satisfies the same properties in~$\R^N\times[0,1-s_1]$. In particular, the functions
$$(t,x)\mapsto h_n(t,x):=\frac{g(x+x_n,v(t_n,x_n)\,v_n(t,x))}{v(t_n,x_n)}$$
are bounded locally in $\R\times\R^N$ and $\|h_n-g_s(x+x_n,0)\,v_n\|_{L^{\infty}(K)}\to0$ as $n\to+\infty$ for any compact set $K\subset\R\times\R^N$, from the mean value theorem. From standard parabolic estimates and Sobolev estimates, it follows that the functions $v_n$ are bounded locally in $W^{1,2,p}_{t,x}(\R\times\R^N)$ and are therefore bounded locally in $C^{0,\alpha}(\R\times\R^N)$. It is then straightforward to check that the functions $h_n$ are actually bounded locally in $C^{0,\alpha}(\R\times\R^N)$. Notice also that, up to extraction of a subsequence, the functions $g_s(\cdot+x_n,0)$ converge locally uniformly in $\R^N$ to a function $a\in C^{0,\alpha}(\R^N)$ and that the matrix fields $A(\cdot+x_n)$ converge locally uniformly in $\R^N$ to a uniformly definite positive symmetric matrix field $A_{\infty}\in C^{1,\alpha}(\R^N)$. As a consequence, again by standard parabolic estimates, the functions $v_n$ converge, up to extraction of a subsequence, locally uniformly in $C^{1,2}_{t,x}(\R\times\R^N)$, to a nonnegative classical solution $v_{\infty}$ of
\be\label{eqvinfty}
(v_{\infty})_t=\hbox{div}(A_{\infty}(x)\nabla v_{\infty})+a(x)\,v_{\infty},\ \ t\in\R,\ x\in\R^N.
\ee\par
On the other hand, $v_n(0,0)=1$, whence $v_{\infty}(0,0)=1$ and $v_{\infty}>0$ in $\R\times\R^N$ from the strong maximum principle. Hence, the functions $(v_n)_t/v_n$ and $\nabla v_n/v_n$ converge locally uniformly in~$\R\times\R^N$ to $(v_{\infty})_t/v_{\infty}$ and $\nabla v_{\infty}/v_{\infty}$. In particular,
$$\frac{(v_{\infty})_t(0,0)}{v_{\infty}(0,0)}=(v_{\infty})_t(0,0)=M.$$
Moreover, since the fields $v_t/t$ and $\nabla v/v$ are bounded in $[t_0,+\infty)\times\R^N$ together with $\lim_{n\to+\infty}t_n=+\infty$ and $\lim_{n\to+\infty}v(t_n,x_n)=0$, it follows that the fields $(v_{\infty})_t/v_{\infty}$ and $\nabla v_{\infty}/v_{\infty}$ are bounded in $\R\times\R^N$ and that $v(t+t_n,x+x_n)\to0$ as $n\to+\infty$ locally uniformly in $\R\times\R^N$. Hence, owing to the definition of $M$ in~\eqref{defM}, one infers that
$$\frac{(v_{\infty})_t(t,x)}{v_{\infty}(t,x)}\le M\ \hbox{ for all }(t,x)\in\R\times\R^N.$$\par
Denote $z=(v_{\infty})_t/v_{\infty}$. Since the coefficients of~\eqref{eqvinfty} do not depend on $t$, it follows from standard parabolic estimates and the differentiation of~\eqref{eqvinfty} with respect to $t$ that $z$ is a classical solution of
\be\label{eqz}
z_t=\hbox{div}(A_{\infty}\nabla z)+2\frac{\nabla v_{\infty}}{v_{\infty}}\cdot A_{\infty}\nabla z
\ee
in $\R\times\R^N$. Furthermore, $z$ and $\nabla v_{\infty}/v_{\infty}$ are bounded in $\R\times\R^N$ and $z\le M$ in $\R\times\R^N$ with~$z(0,0)=M$. The strong parabolic maximum principle then implies that $z=M$ in~$(-\infty,0]\times\R^N$, and hence in $\R\times\R^N$ by uniqueness of the Cauchy problem associated with~\eqref{eqz}. In other words, $(v_{\infty})_t/v_{\infty}=M$ in $\R\times\R^N$. In particular, since $v_{\infty}(0,0)=1$, one gets that~$v_{\infty}(\tau,0)=e^{M\tau}$ (we recall that $\tau=\tau_*+1$).\par
Lastly, by Lemma~\ref{lem3} applied with $\tau=\tau_*+1$ and $a=1/2$, there are $T_0>0$ and $\delta>0$ such that
$$1-u(t+\tau,x)\le(1-\delta)\,(1-u(t,x))\ \hbox{ for all }(t,x)\in[T_0,+\infty)\times\R^N\hbox{ with }u(t,x)\ge\frac12.$$
Thus, for any given $(t,x)\in\R\times\R^N$, since $t+t_n\ge T_0$ and
$$u(t+t_n,x+x_n)=1-v(t+t_n,x+x_n)\ge\frac12$$
for all $n$ large enough, one infers that
$$1-u(t+t_n+\tau,x+x_n)\le(1-\delta)\,(1-u(t+t_n,x+x_n)),$$
whence $v_n(t+\tau,x)\le(1-\delta)\,v_n(t,x)$ for all $n$ large enough. Thus, $v_{\infty}(t+\tau,x)\le(1-\delta)\,v_{\infty}(t,x)$ for all $(t,x)\in\R\times\R^N$. Consequently, $e^{M\tau}=v_{\infty}(\tau,0)\le(1-\delta)v_{\infty}(0,0)=1-\delta<1$ and $M<0$.\par
As a conclusion, owing to the definition of $M$ in~\eqref{defM} and since $v=1-u$, the conclusion of Lemma~\ref{lem4} follows.\hfill$\Box$


\SE{$\tau$-monotonicity in time when $u(t,x)$ is close to $0$}\label{sec6}

In this section, for any arbitrary $\tau>0$, we show the $\tau$-monotonicity in time at large time in the region where $u(t,x)$ is close to $0$. We shall use in particular the assumptions~\eqref{oscA} and~\eqref{oscf} on asymptotic homogeneity of the coefficients $A$ and $f_u(\cdot,0)$. The key step will be the following proposition, which is of independent interest.

\begin{pro}\label{pro1}
Let $\underline{\nu}$ and $\overline{\nu}$ be two fixed positive real numbers such that $0<\underline{\nu}\le\overline{\nu}$ and let~$\sigma\in(0,1)$ be fixed. Then, there exist $\tau>0$ and $\eta>0$ such that, for every $C^1(\R^N)$ symmetric matrix field $a=(a_{ij})_{1\le i,j\le N}:\R^N\to\mathbb{S}_N(\R)$ with $\underline{\nu}I\le a\le\overline{\nu}I$ and $|\nabla a|\le\eta$ in $\R^N$ $($where~$|\nabla a(x)|=\max_{1\le i,j\le N}|\nabla a_{ij}(x)|)$, the fundamental solution $p(t,x;y)$ of
\be\label{defp}\left\{\baa{rcl}
p_t(t,x;y) & = & {\rm{div}}(a(x)\nabla p(t,x;y)),\ \ t>0,\ x\in\R^N,\vspace{3pt}\\
p(0,\cdot;y) & = & \delta_y\eaa\right.
\ee
satisfies
\be\label{ptau}
p(\tau+1,x;0)\ge\sigma\,p(\tau,x;0)\ \hbox{ for all }x\in\R^N.
\ee
\end{pro}

Let us postpone the proof of this proposition to Section~\ref{secpro1}. We continue in this section the proof of Theorem~\ref{th1} for the solution $u$ of~\eqref{eq}. The main result proved in this section is the following lemma.

\begin{lem}\label{lem5}
Let $\theta$ and $\theta'$ be any two real numbers such that $0<\theta\le\theta'$. Then there exist $T_0>0$ and $\epsilon>0$ such that
$$\forall\,\tau\in[\theta,\theta'],\ \forall\,(t,x)\in[T_0,+\infty)\times\R^N,\ \ \big(u(t,x)\le\epsilon\big)\Longrightarrow\big(u(t+\tau,x)>u(t,x)\big).$$
\end{lem}

\noindent{\bf{Proof.}} Let us argue by way of contradiction. So, assume that there exist a sequence $(\tau_n)_{n\in\N}$ in~$[\theta,\theta']$, a sequence $(t_n)_{n\in\N}$ of positive real numbers and a sequence $(x_n)_{n\in\N}$ of points in $\R^N$ such that
\be\label{tntaun}
t_n\mathop{\longrightarrow}_{n\to+\infty}+\infty,\ u(t_n,x_n)\mathop{\longrightarrow}_{n\to+\infty}0\hbox{ and }u(t_n+\tau_n,x_n)\le u(t_n,x_n)\hbox{ for all }n\in\N.
\ee
Up to extraction of a subsequence, one can assume without loss of generality that
$$\tau_n\to\tau_{\infty}\in[\theta,\theta']\subset(0,+\infty)\ \hbox{ as }n\to+\infty.$$
Notice that~\eqref{defc} and~\eqref{tntaun} yield $\lim_{n\to+\infty}|x_n|=+\infty$ and even $\liminf_{n\to+\infty}|x_n|/t_n\ge c>0$. Therefore, up to extraction of another subsequence, the matrix fields $x\mapsto A(x+x_n)$ converge in~$C^1_{loc}(\R^N)$ to a definite positive symmetric matrix field $A_{\infty}$, which turns out to be a constant matrix due to~\eqref{oscA}. Observe also that~\eqref{tntaun} and the nonnegativity of $u$ imply that~$u(t_n+\tau_n,x_n)\to0$ as $n\to+\infty$.\par
Define the functions
$$u_n(t,x)=\frac{u(t+t_n+\tau_n,x+x_n)}{u(t_n+\tau_n,x_n)},\ \ (t,x)\in(-t_n-\tau_n,+\infty)\times\R^N.$$
Each function $u_n$ obeys
$$(u_n)_t(t,x)=\hbox{div}(A(x+x_n)\nabla u_n(t,x))+\frac{f(x+x_n,u(t_n+\tau_n,x_n)\,u_n(t,x))}{u(t_n+\tau_n,x_n)}.$$
For each compact set $K\subset(-\infty,0)\times\R^N$, there is $n_K\in\N$ such that $K\subset(-t_n-\tau_n+1,+\infty)\times\R^N$ for all $n\!\ge\!n_K$ and it follows from Harnack inequality applied to $u$ that $\sup_{n\ge n_K}\!\|u_n\|_{L^{\infty}(K)}\!<\!+\!\infty$. Remember that $f(\cdot,0)=0$ in $\R^N$, that the function $(x,s)\mapsto f(x,s)$ is Lipschitz continuous with respect to $s$ uniformly in $x\in\R^N$, of class $C^1$ with respect to $s$ in $\R^N\times[0,s_0]$ with~$s_0\in(0,1)$, and that $f_s$ is uniformly continuous in $\R^N\times[0,s_0]$ and of class $C^{0,\alpha}$ with respect to~$x$ uniformly in~$s\in[0,s_0]$. Furthermore, up to extraction of another subsequence, the functions~$x\mapsto f_s(x+x_n,0)$ converge locally uniformly in $\R^N$ to a function $r\in C^{0,\alpha}(\R^N)$, which is actually a constant such that~$r\ge\mu>0$ by~\eqref{hypf} and~\eqref{oscf}. Therefore, as we did in the proof of Lemma~\ref{lem4} for the functions~$v_n$ satisfying~\eqref{eqvn}, we get that, up to extraction of a subsequence, the positive functions~$u_n$ converge locally uniformly in $C^{1,2}_{t,x}((-\infty,0)\times\R^N)$ to a nonnegative solution $u_{\infty}$ of
$$(u_{\infty})_t=\hbox{div}(A_{\infty}\nabla u_{\infty})+r\,u_{\infty}\ \hbox{ in }(-\infty,0)\times\R^N.$$
Since $u_n(-\tau_n,0)\ge1$ by~\eqref{tntaun} and $\tau_n\to\tau_{\infty}>0$ as $n\to+\infty$, we get that $u_{\infty}(-\tau_{\infty},0)\ge1$, whence $u_{\infty}>0$ in $(-\infty,0)\times\R^N$ from the strong maximum principle and the uniqueness of the Cauchy problem associated with that equation. Since $A_{\infty}$ is a constant symmetric definite positive matrix, there is an invertible matrix $M$ such that the function $\tilde{u}_{\infty}$ defined in $(-\infty,0)\times\R^N$ by $\tilde{u}_{\infty}(t,x)=u_{\infty}(t,Mx)$ satisfies $(\tilde{u}_{\infty})_t=\Delta\tilde{u}_{\infty}+r\,\tilde{u}_{\infty}$ in $(-\infty,0)\times\R^N$. In other words, the function $(t,x)\mapsto e^{-rt}\tilde{u}_{\infty}(t,x)$ is a positive solution of the heat equation in $(-\infty,0)\times\R^N$. Thus, by~\cite{w}, it is nondecreasing with respect to $t$ in $(-\infty,0)\times\R^N$. Therefore, the function $(t,x)\mapsto e^{-rt}u_{\infty}(t,x)$ is nondecreasing with respect to $t$ in $(-\infty,0)\times\R^N$.\par
Remember now that $\mu>0$ and that $\theta>0$ is given in the statement of Lemma~\ref{lem5}. Fix a real number $\sigma\in(0,1)$ close enough to $1$ so that
$$\sigma\,e^{\,\mu\,\theta/2}>1,$$
and let
$$\tau>0\hbox{ and }\eta>0$$
be as in Proposition~\ref{pro1} applied with this real number $\sigma$ and with $\underline{\nu}=\nu^{-1}$ and $\overline{\nu}=\nu$ (remember that~$\nu^{-1}I\le A\le\nu I$ in $\R^N$ with $\nu\ge1$). Finally, let $L>0$ be such that~\eqref{defL} holds and let us fix~$\epsilon>0$ small enough so that
\be\label{defepsilon}
0<\epsilon\le\frac{\theta}{4},\ \ \sqrt{\epsilon}\,|\nabla A|\le\eta\hbox{ in }\R^N\ \hbox{ and }\ e^{-\epsilon L\tau}\,\sigma\,e^{\,\mu\,\theta/2}>1.
\ee\par
Let us finally complete the proof of Lemma~\ref{lem5} by reaching a contradiction. Since the function~$t\mapsto e^{-rt}u_{\infty}(t,x)$ is nondecreasing in $(-\infty,0)$ for each given $x\in\R^N$ and since~$0<\epsilon<2\epsilon<\theta\le\tau_{\infty}$, one has $e^{r\epsilon}u_{\infty}(-\epsilon,0)\ge e^{r\tau_{\infty}}u_{\infty}(-\tau_{\infty},0)$. But $u_n\to u_{\infty}>0$ locally uniformly in $(-\infty,0)\times\R^N$ and $\tau_n\to\tau_{\infty}$ as $n\to+\infty$. As a consequence,
$$e^{r\epsilon}u_n(-\epsilon,0)\ge e^{r(\tau_n-\epsilon)}u_n(-\tau_n,0)\ \hbox{ for all }n\hbox{ large enough},$$
whence
$$u(-\epsilon+t_n+\tau_n,x_n)\ge e^{r(\tau_n-2\epsilon)}u(t_n,x_n)\ge e^{\,\mu\,\theta/2}u(t_n,x_n)\ge e^{\,\mu\,\theta/2}u(t_n+\tau_n,x_n)$$
for all $n$ large enough, since $r\ge\mu>0$, $\tau_n-2\epsilon\ge\theta-2\epsilon\ge\theta/2>0$ and $u(t_n,x_n)\ge u(t_n+\tau_n,x_n)>0$ by~\eqref{tntaun}. Lastly, consider the parabolically rescaled functions
$$v_n(t,x)=u(\epsilon\,t+t_n+\tau_n,\sqrt{\epsilon}\,x+x_n),\ \ (t,x)\in(-\epsilon^{-1}(t_n+\tau_n),+\infty)\times\R^N$$
and observe that
\be\label{vn10}
v_n(-1,0)\ge e^{\,\mu\,\theta/2}v_n(0,0)\ \hbox{ for all }n\hbox{ large enough}.
\ee
Furthermore, the functions $v_n$ obey
\be\label{eqvnbis}
(v_n)_t=\hbox{div}(A_n(x)\nabla v_n)+\epsilon\,f(\sqrt{\epsilon}\,x+x_n,v_n)\ \hbox{ in }(-\epsilon^{-1}(t_n+\tau_n),+\infty)\times\R^N,
\ee
with $A_n(x)=A(\sqrt{\epsilon}\,x+x_n)$, and they are positive in $(-\epsilon^{-1}(t_n+\tau_n),+\infty)\times\R^N$. For each $n\in\N$ and $y\in\R^N$, call $p_n$ and $p_{n,y}$ the fundamental solutions of~\eqref{defp} with diffusion matrix fields $a=A_n$ and $a=A_n(\cdot+y)$, respectively. Remember that $\tau>0$ is given above from Proposition~\ref{pro1} and choose~$n\in\N$ large enough so that $-\epsilon^{-1}(t_n+\tau_n)<-\tau-1$ and~\eqref{vn10} holds. Since the function $f$ is such that~$0\le f(x,s)\le Ls$ for all $(x,s)\in\R^N\times[0,1]$, it follows from~\eqref{eqvnbis} that
\be\label{vn00}
v_n(0,0)\ge\int_{\R^N}p_n(\tau+1,0;y)\,v_n(-\tau+1,y)\,dy=\int_{\R^N}p_{n,y}(\tau+1,-y;0)\,v_n(-\tau+1,y)\,dy
\ee
and
\be\label{vn-10}
v_n(-1,0)\le e^{\epsilon L\tau}\int_{\R^N}p_n(\tau,0;y)\,v_n(-\tau+1,y)\,dy=e^{\epsilon L\tau}\int_{\R^N}p_{n,y}(\tau,-y;0)\,v_n(-\tau+1,y)\,dy.
\ee
On the other hand, for every $y\in\R^N$, the matrix field $A_{n,y}:=A_n(\cdot+y)$ is of class $C^1(\R^N)$ it satisfies $\nu^{-1}I\le A_{n,y}\le\nu I$ in $\R^N$ and~$|\nabla A_{n,y}(x)|=\sqrt{\epsilon}\,|\nabla A(\sqrt{\epsilon}\,(x+y)+x_n)|\le\eta$ for all $x\in\R^N$ by~\eqref{defepsilon}, where $\eta>0$ is given above from Proposition~\ref{pro1}. It follows then from the conclusion of Proposition~\ref{pro1} that
$$p_{n,y}(\tau+1,-y;0)\ge\sigma\,p_{n,y}(\tau,-y;0)\ \hbox{ for all }y\in\R^N,$$
whence $v_n(0,0)\ge e^{-\epsilon L\tau}\sigma\,v_n(-1,0)$ by~\eqref{vn00} and~\eqref{vn-10}, and finally $v_n(0,0)\ge e^{-\epsilon L\tau}\sigma\,e^{\,\mu\,\theta/2}v_n(0,0)$ by~\eqref{vn10}. The posi\-tivity of $v_n(0,0)=u(t_n+\tau_n,x_n)$ (since $t_n+\tau_n>0$) contradicts the last property in~\eqref{defepsilon}. The proof of Lemma~\ref{lem5} is thereby complete.\hfill$\Box$


\SE{Conclusion of the proof of Theorem~\ref{th1}}\label{sec7}

Remember the definition of $\tau_*$ in~\eqref{deftau*} and remember that $0\le\tau^*<+\infty$. Before completing the proof of Theorem~\ref{th1}, we first show that $\tau^*=0$.

\begin{lem}\label{lem6}
There holds
$$\tau_*=0.$$
\end{lem}

\noindent{\bf{Proof.}} Assume by contradiction that $\tau_*>0$. It follows from the definition~\eqref{deftau*} of $\tau_*$ that there exist two sequences $(\tau_n)_{n\in\N}$ and $(t_n)_{n\in\N}$ of positive real numbers and a sequence $(x_n)_{n\in\N}$ of points in $\R^N$ such that
\be\label{tntaunu}
\lim_{n\to+\infty}t_n=+\infty,\ \ \liminf_{n\to+\infty}\tau_n\ge\tau_*\ \hbox{ and }\ u(t_n+\tau_n,x_n)<u(t_n,x_n)\hbox{ for all }n\in\N.
\ee
Actually, it turns out that $\tau_n\to\tau_*$ as $n\to+\infty$, otherwise the assumption $\limsup_{n\to+\infty}\tau_n>\tau_*$ together with $\lim_{n\to+\infty}t_n=+\infty$ would contradict the definition of $\tau_*$. The real numbers $u(t_n,x_n)$ take values in $(0,1)$. Thus, up to extraction of a subsequence, three cases can occur.\par
{\it Case 1: $u(t_n,x_n)\to m\in(0,1)$ as $n\to+\infty$.} There are then two real numbers $a$ and $b$ such that~$0<a\le u(t_n,x_n)\le b<1$ for all $n\in\N$. Since $t_n\to+\infty$ and $\tau_n\to\tau_*>0$ as $n\to+\infty$, Corollary~\ref{cor2} yields the existence of $\delta>0$ such that $u(t_n+\tau_n,x_n)\ge(1+\delta)\,u(t_n,x_n)$ for all $n$ large enough. This is impossible by~\eqref{tntaunu}.\par
{\it Case 2: $u(t_n,x_n)\to1$ as $n\to+\infty$.} Therefore, $u(t_n,x_n)\ge b$ and $t_n\ge\tilde{T}$ for all $n$ large enough, where $b\in(0,1)$ and $\tilde{T}$ are given as in Lemma~\ref{lem4}. One infers then from Lemma~\ref{lem4} that, for all $n$ large enough, $u_t(t_n,x_n)>0$ and thus even that $u_t(t,x_n)>0$ for all $t\ge t_n$. In particular,~$u(t_n+\tau_n,x_n)>u(t_n,x_n)$ for all $n$ large enough, contradicting~\eqref{tntaunu}. Thus, Case~2 is ruled out too.\par
{\it Case 3: $u(t_n,x_n)\to0$ as $n\to+\infty$.} Since $\tau_n\to\tau_*>0$ and each $\tau_n$ is positive, there are two real numbers $0<\theta\le\theta'$ such that $\theta\le\tau_n\le\theta'$ for all $n\in\N$. Let then $T_0>0$ and~$\epsilon>0$ be as in Lemma~\ref{lem5}. For all $n$ large enough, one has $t_n\ge T_0$ and $u(t_n,x_n)\le\epsilon$, whence~$u(t_n+\tau_n,x_n)>u(t_n,x_n)$. This contradicts~\eqref{tntaunu}.\par
As a conclusion, all three cases are impossible. Finally, $\tau_*=0$ and the proof of Lemma~\ref{lem6} is complete.\hfill$\Box$\break

Based on the previous lemma, we can now complete the proof of Theorem~\ref{th1}.\hfill\break

\noindent{\bf{Proof of Theorem~\ref{th1}.}} We first notice that, for any $t>0$,
\be\label{ut0}
u_t(t,x)\to0\ \hbox{ as }|x|\to+\infty.
\ee
Indeed, on the one hand, as in the proof of Lemma~\ref{lem2}, for any sequence $(x_n)_{n\in\N}$ in $\R^N$ with~$\lim_{n\to+\infty}|x_n|=+\infty$, the functions $u_n:(t,x)\mapsto u(t,x+x_n)$ converge locally uniformly in $C^{1,2}_{t,x}((0,+\infty)\times\R^N)$, up to extraction of a subsequence, to a solution $0\le u_{\infty}\le 1$ of an equation of the type
$$(u_{\infty})_t=\hbox{div}(A_{\infty}\nabla u_{\infty})+f_{\infty}(x,u_{\infty})\ \hbox{ in }(0,+\infty)\times\R^N$$
for some constant symmetric definite positive matrix $A_{\infty}$ and for some function $f_{\infty}$ satis\-fying the same properties as $f$. On the other hand, as already emphasized from the proof of Lemma~\ref{lem1},~$u(t,x)\to0$ as $|x|\to+\infty$, for every $t>0$. Therefore, $u_{\infty}=0$ in $(0,+\infty)\times\R^N$. Hence, by uniqueness the whole sequence $(u_n)_{n\in\N}$ converges to $0$ locally uniformly in $C^{1,2}_{t,x}((0,+\infty)\times\R^N)$ and~$u_t(t,x+x_n)\to0$ as $n\to+\infty$ for every $(t,x)\in(0,+\infty)\times\R^N$. As a consequence,~\eqref{ut0} holds.\par
In particular, it follows that $\inf_{\R^N}u_t(t,\cdot)\le0$ for all $t>0$. We now prove~\eqref{ut}, and then~\eqref{Teps}.\footnote{We could also view~\eqref{ut} as a consequence of~\eqref{Teps} by observing that $u_t(t,x)\to0$ as $u(t,x)\to0$. But since the proof of~\eqref{ut} is easy even without~\eqref{Teps}, we choose to carry it out before~\eqref{Teps}.} Assume now by contradiction that~$\inf_{\R^N}u_t(t,\cdot)\not\to0$ as $t\to+\infty$. Since $u_t$ is bounded in~$(1,+\infty)\times\R^N$, it follows then that there are a sequence $(t_n)_{n\in\N}$ of positive real numbers and a sequence $(x_n)_{n\in\N}$ of points in $\R^N$ such that~$t_n\to+\infty$ and $\limsup_{n\to+\infty}u_t(t_n,x_n)\in(-\infty,0)$. As done in the previous paragraph or in the proof of Lemma~\ref{lem2}, the functions
$$(t,x)\mapsto u(t+t_n,x+x_n)$$
converge, up to extraction of a subsequence, locally uniformly in $C^{1,2}_{t,x}(\R\times\R^N)$ to a solution~$0\le u_{\infty}\le1$ of an equation of the type~\eqref{uinfty} with $(u_{\infty})_t(0,0)<0$. However, for any $\tau>0$, it follows from the definition~\eqref{deftau*} of $\tau_*$ together with $\lim_{n\to+\infty}t_n=+\infty$ and Lemma~\ref{lem6} ($\tau_*=0$) that, for any given~$(t,x)\in\R^N$, there holds
$$u(t+\tau+t_n,x+x_n)\ge u(t+t_n,x+x_n)\ \hbox{ for all }n\hbox{ large enough},$$
whence $u_{\infty}(t+\tau,x)\ge u_{\infty}(t,x)$. Therefore, since $\tau>0$ and $(t,x)\in\R\times\R^N$ are arbitrary, one gets that $(u_{\infty})_t\ge0$ in $\R\times\R^N$, which yields a contradiction. In other words, one has shown that~$\inf_{\R^N}u_t(t,\cdot)\to0$ as $t\to+\infty$, that is,~\eqref{ut}.\par
Finally, let $\epsilon\in(0,1)$ be given and let us show~\eqref{Teps}, that is, the existence of $T_{\epsilon}>0$ such that~$u_t(t,x)>0$ for every $(t,x)\in[T_{\epsilon},+\infty)\times\R^N$ with $u(t,x)\ge\epsilon$. Assume not. Then there are a sequence~$(t_n)_{n\in\N}$ of positive real numbers and a sequence $(x_n)_{n\in\N}$ of points in $\R^N$ such that~$t_n\to+\infty$ as~$n\to+\infty$, while $u(t_n,x_n)\ge\epsilon$ and $u_t(t_n,x_n)\le0$ for all $n\in\N$. It follows then from Lemma~\ref{lem4} that~$1$ is not a limiting value of the sequence $(u(t_n,x_n))_{n\in\N}$. Therefore, up to extraction of a sequence, one can assume without loss of generality that
$$u(t_n,x_n)\to m\in(0,1)\ \hbox{ as }n\to+\infty.$$
As done in the previous paragraph, one infers that the functions $(t,x)\mapsto u(t+t_n,x+x_n)$ converge, up to extraction of a subsequence, locally uniformly in $C^{1,2}_{t,x}(\R\times\R^N)$ to a solution $0\le u_{\infty}\le1$ of an equation of the type
$$(u_{\infty})_t=\hbox{div}(A_{\infty}\nabla u_{\infty})+f_{\infty}(x,u_{\infty}),\ \ t\in\R,\ x\in\R^N$$
with $(u_{\infty})_t(t,x)\ge0$ for all $(t,x)\in\R\times\R^N$, whereas $u_{\infty}(0,0)=m\in(0,1)$ and $(u_{\infty})_t(0,0)\le0$. It follows from the strong maximum principle applied to the function $(u_{\infty})_t$ that $(u_{\infty})_t=0$ in~$(-\infty,0]\times\R^N$ and then in $\R\times\R^N$. Furthermore, the strong maximum principle applied to~$u_{\infty}$ also implies that $0<u_{\infty}<1$ in $\R\times\R^N$. As recalled in Section~\ref{intro}, one infers then from the properties shared by $f_{\infty}$ with $f$ that $u_{\infty}(t,x)\to1$ as $t\to+\infty$ for every $x\in\R^N$. This leads to a contradiction, since $(u_{\infty})_t=0$ in $\R\times\R^N$ and $u_{\infty}(0,0)=m<1$. Therefore,~\eqref{Teps} is shown and the proof of Theorem~\ref{th1} is thereby complete.\hfill$\Box$
 

\SE{Proof of Proposition~\ref{pro1}}\label{secpro1}

Let $0<\underline{\nu}\le\overline{\nu}$ and $0<\sigma<1$ be fixed. When $a=D\,I$ with a real number $D\in[\underline{\nu},\overline{\nu}]$, then the conclusion~\eqref{ptau} holds immediately for $\tau>0$ large enough, from the explicit formula
$$p(t,x;0)=\frac{e^{-|x|^2/(4Dt)}}{(4\pi t)^{N/2}}.$$
In the general case where $a$ may not be constant, we will get the estimates~\eqref{ptau} by using uniform Gaussian estimates for large $x$ and small $t$, and by approximating locally, when $a$ is nearly locally constant, the solutions $p$ of~\eqref{defp} by explicit fundamental solutions of parabolic equations with constant coefficients.\par
More precisely, choose first any real number $\tau>0$ large enough so that
\be\label{tausigma}
\Big(\frac{\tau}{\tau+1}\Big)^{N/2}>\sigma.
\ee
We shall prove the conclusion of Proposition~\ref{pro1} is fulfilled with any such real number $\tau$.\par
First of all, it follows from the Gaussian upper and lower bounds of the fundamental solutions of~\eqref{defp} in~\cite{n} that there exist a constant $K\ge 1$ such that, for every $L^{\infty}(\R^N)$ matrix field~$a:\R^N\to\mathbb{S}_N(\R)$ with $\underline{\nu}I\le a\le\overline{\nu}I$ a.e. in $\R^N$, the fundamental solution $p$ of~\eqref{defp} satisfies
\be\label{pbounds}
\frac{e^{-K\,|x|^2/t}}{K\,t^{N/2}}\le p(t,x;0)\le\frac{K\,e^{-|x|^2/(Kt)}}{t^{N/2}}\ \hbox{ for all }(t,x)\in(0,+\infty)\times\R^N.
\ee\par
Therefore, we claim that there exists $\tau_0\in(0,\tau)$ small enough such that, for every $L^{\infty}(\R^N)$ matrix field $a:\R^N\to\mathbb{S}_N(\R)$ with $\underline{\nu}I\le a\le\overline{\nu}I$ a.e. in $\R^N$, there holds
\be\label{ptau0}
p(\tau_0+1,x;0)\ge\sigma\,p(\tau_0,x;0)\ \hbox{ for all }x\in\R^N\hbox{ with }|x|\ge1.
\ee
Indeed, otherwise, by picking any sequence $(\tau_n)_{n\in\N}$ in $(0,\tau)$ such that $\lim_{n\to+\infty}\tau_n=0$, there would exist a sequence $(a_n)_{n\in\N}$ of bounded symmetric matrix fields with $\underline{\nu}I\le a_n\le\overline{\nu}I$ a.e. in $\R^N$ and a sequence $(x_n)_{n\in\N}$ of points in $\R^N$ such that
$$|x_n|\ge1\ \hbox{ and }\ p_n(\tau_n+1,x_n;0)<\sigma\,p_n(\tau_n,x_n;0)\ \hbox{ for all }n\in\N,$$
where $p_n$ denotes the fundamental solution of~\eqref{defp} with $a=a_n$. It would then follow from~\eqref{pbounds} that
$$\frac{e^{-K|x_n|^2/(\tau_n+1)}}{K\,(\tau_n+1)^{N/2}}<\frac{\sigma\,K\,e^{-|x_n|^2/(K\tau_n)}}{\tau_n^{N/2}}$$
for all $n\in\N$, that is,
$$e^{|x_n|^2(1/(K\tau_n)-K/(\tau_n+1))}<\sigma\,K^2\,\Big(\frac{\tau_n+1}{\tau_n}\Big)^{N/2}.$$
But $|x_n|\ge1$ for all $n\in\N$ and $1/(K\tau_n)-K/(\tau_n+1)\ge1/(2K\tau_n)$ for all $n$ large enough, since~$\tau_n\to0^+$ as $n\to+\infty$. Therefore, $e^{1/(2K\tau_n)}<\sigma\,K^2\,(1+1/\tau_n)^{N/2}$ for all $n$ large enough, which leads to a contradiction by passing to the limit as $n\to+\infty$. As a consequence, there is $\tau_0\in(0,\tau)$ such that~\eqref{ptau0} holds for every $L^{\infty}(\R^N)$ matrix field~$a:\R^N\to\mathbb{S}_N(\R)$ with $\underline{\nu}I\le a\le\overline{\nu}I$ a.e. in~$\R^N$.\par
Now, we fix a real number $R\ge1$ large enough such that
\be\label{defR}
\sigma\,\Big(1+\frac{1}{\tau_0}\Big)^{N/2}e^{-R^2/(4\overline{\nu}\tau(\tau+1))}<1.
\ee
Given this choice of $R\ge1$, we claim that there exists a real number $\eta>0$ such that, for every~$a\in C^1(\R^N;\mathbb{S}_N(\R))$ with $\underline{\nu}I\le a\le\overline{\nu}I$ and $|\nabla a|\le\eta$ in $\R^N$, the fundamental solution~$p(t,x;y)$ of~\eqref{defp} satisfies
\be\label{boundaries}\left\{\baa{rl}
p(\tau+1,x;0)\ge\sigma\,p(\tau,x;0) & \hbox{for all }|x|\le R,\vspace{3pt}\\
p(t+1,x;0)\ge\sigma\,p(t,x;0) & \hbox{for all }t\in[\tau_0,\tau]\hbox{ and }|x|=R.\eaa\right.
\ee
Indeed, otherwise, there exist a sequence $(a_n)_{n\in\N}$ in $C^1(\R^N;\mathbb{S}_N(\R))$ with $\underline{\nu}I\le a_n\le\overline{\nu}I$ in $\R^N$ and~$\lim_{n\to+\infty}\|\,|\nabla a_n|\,\|_{L^{\infty}(\R^N)}=0$, as well as a sequence of points $(t_n,x_n)_{n\in\N}$ in $(0,+\infty)\times\R^N$ such that
\be\label{tnxn}
(t_n,x_n)\,\in\,\{\tau\}\!\times\!\overline{B(0,R)}\,\cup\,[\tau_0,\tau]\!\times\!\partial B(0,R)\ \hbox{ and }\ p_n(t_n+1,x_n;0)<\sigma\,p_n(t_n,x_n;0)
\ee
for all $n\in\N$, where $p_n$ denotes the fundamental solution of~\eqref{defp} with $a=a_n$. Up to extraction of a subsequence, the matrix fields $a_n$ converge locally uniformly in $\R^N$ to a constant symmetric definite positive matrix $a_{\infty}$ such that $\underline{\nu}I\le a_{\infty}\le\overline{\nu}I$. Furthermore, the functions $(p_n)_{n\in\N}$ are bounded locally in $(0,+\infty)\times\R^N$ from the bounds~\eqref{pbounds}. From standard parabolic estimates, the functions $p_n(\cdot,\cdot;0)$ converge then locally uniformly in $(0,+\infty)\times\R^N$ to a classical solution~$p_{\infty}$ of
$$(p_{\infty})_t=\hbox{div}(a_{\infty}\nabla p_{\infty})\ \hbox{ in }(0,+\infty)\times\R^N$$
such that
\be\label{pinfty0}
K^{-1}t^{-N/2}e^{-K|x|^2/t}\le p_{\infty}(t,x)\le K\,t^{-N/2}e^{-|x|^2/(Kt)}\ \hbox{ for all }(t,x)\in(0,+\infty)\times\R^N.
\ee
Moreover, it follows from~\eqref{tnxn} that there exists a point $(t_{\infty},x_{\infty})$ such that
\be\label{tinfty}
(t_{\infty},x_{\infty})\,\in\,\{\tau\}\!\times\!\overline{B(0,R)}\,\cup\,[\tau_0,\tau]\!\times\!\partial B(0,R)\ \hbox{ and }\ p_{\infty}(t_{\infty}+1,x_{\infty})\le\sigma\,p_{\infty}(t_{\infty},x_{\infty}).
\ee
From the aforementioned Gaussian estimates, the function $p_{\infty}$ is therefore positive in $(0,+\infty)\times\R^N$ and since $a_{\infty}\in\mathbb{S}_N(\R^N)$ satisfies $\underline{\nu}I\le a_{\infty}\le\overline{\nu}I$, there is an orthogonal linear map $M:x\mapsto y$ such that the function $(t,y)\mapsto q(t,y)=p_{\infty}(t,M^{-1}y)=p_{\infty}(t,x)$ is a positive solution of
$$q_t=\sum_{1\le i\le N}\lambda_i\frac{\partial^2q}{\partial y_i^2}\ \hbox{ in }(0,+\infty)\times\R^N$$
for some real numbers $\lambda_1,\ldots,\lambda_N\in[\underline{\nu},\overline{\nu}]$. Hence, there is a nonnegative Radon measure $\lambda$ such that
$$p_{\infty}(t,x)=\frac{1}{(4\pi t)^{N/2}}\int_{\R^N}e^{-\sum_{1\le i\le N}|y_i/\sqrt{\lambda_i}-z_i|^2/(4t)}\,d\lambda(z)\ \hbox{ for all }(t,x)\in(0,+\infty)\times\R^N.$$
Since, by~\eqref{pinfty0}, $\int_{B(x_0,r)}p_{\infty}(t,x)\,dx\to0^+$ as $t\to0^+$ for every $x_0\in\R^N$ and $r>0$ such that~$0\not\in\overline{B(x_0,r)}$, it follows easily that $\lambda$ is supported on the singleton~$\{0\}$ and that there is~$\rho>0$ such that
\be\label{pinfty}
p_{\infty}(t,x)=\frac{\rho\,e^{-\sum_{1\le i\le N}|y_i|^2/(4\lambda_it)}}{(4\pi t)^{N/2}}\ \hbox{ for all }(t,x)\in(0,+\infty)\times\R^N.
\ee
Remember now~\eqref{tinfty}. On the one hand, if $t_{\infty}=\tau$ (and $|x_{\infty}|\le R$), then~\eqref{tinfty} and~\eqref{pinfty} imply that
$$\frac{\rho\,e^{-\sum_{1\le i\le N}|y_{\infty,i}|^2/(4\lambda_i(\tau+1))}}{(4\pi(\tau+1))^{N/2}}\le\sigma\times\frac{\rho\,e^{-\sum_{1\le i\le N}|y_{\infty,i}|^2/(4\lambda_i\tau)}}{(4\pi\tau)^{N/2}},$$
where $y_{\infty}=M\,x_{\infty}$. Since $\rho>0$, $\lambda_i>0$ for all $1\le i\le N$ and $0<\tau<\tau+1$, one gets that~$\tau^{N/2}\le\sigma(\tau+1)^{N/2}$, which contradicts~\eqref{tausigma}. On the other hand, if $|x_{\infty}|=R$ (and $t_{\infty}\in[\tau_0,\tau]$), then~\eqref{tinfty} and~\eqref{pinfty} yield, with the same notations as above,
$$\frac{\rho\,e^{-\sum_{1\le i\le N}|y_{\infty,i}|^2/(4\lambda_i(t_{\infty}+1))}}{(4\pi(t_{\infty}+1))^{N/2}}\le\sigma\times\frac{\rho\,e^{-\sum_{1\le i\le N}|y_{\infty,i}|^2/(4\lambda_it_{\infty})}}{(4\pi t_{\infty})^{N/2}},$$
whence
$$1\le\sigma\,\Big(1+\frac{1}{t_{\infty}}\Big)^{N/2}e^{-\sum_{1\le i\le N}|y_{\infty,i}|^2/(4\lambda_it_{\infty}(t_{\infty}+1))}.$$
Since $0<\tau_0\le t_{\infty}\le\tau$, $0<\underline{\nu}\le\lambda_i\le\overline{\nu}$ for all $1\le i\le N$ and $|y_{\infty}|=|M\,x_{\infty}|=|x_{\infty}|=R$, one infers that
$$1\le\sigma\,\Big(1+\frac{1}{\tau_0}\Big)^{N/2}e^{-R^2/(4\overline{\nu}\tau(\tau+1))},$$
which contradicts~\eqref{defR}.\par
As a consequence, there is $\eta>0$ such that~\eqref{boundaries} holds for every $a\in C^1(\R^N;\mathbb{S}_N(\R))$ with~$\underline{\nu}I\le a\le\overline{\nu}I$ and $|\nabla a|\le\eta$ in $\R^N$. Consider finally any such matrix field $a$ and let us show that the conclusion~\eqref{ptau} holds, that is $p(\tau+1,\cdot;0)\ge\sigma\,p(\tau,\cdot;0)$ in $\R^N$, where $p$ solves~\eqref{defp}. First of all, it follows from~\eqref{boundaries} that
\be\label{ptau1}
p(\tau+1,\cdot;0)\ge\sigma\,p(\tau,\cdot;0)\ \hbox{ in }\overline{B(0,R)}.
\ee
On the other hand,
$$p(\tau_0+1,\cdot;0)\ge\sigma\,p(\tau_0,\cdot;0)\ \hbox{ in }\R^N\backslash B(0,R)\subset\R^N\backslash B(0,1)$$
by~\eqref{ptau0} and $R\ge 1$. Lastly,
$$p(t+1,\cdot;0)\ge\sigma\,p(t,\cdot;0)\ \hbox{ on }\partial B(0,R)\ \hbox{ for all }t\in[\tau_0,\tau]$$
by~\eqref{boundaries}. Therefore, since $p(\cdot,\cdot;0)$ and $p(\cdot+1,\cdot;0)$ are two positive bounded solutions of the same linear parabolic equation in (at least) $[\tau_0,\tau]\times(\R^N\backslash B(0,R))$, it follows from the parabolic maximum principle that
$$p(\tau+1,\cdot;0)\ge\sigma\,p(\tau,\cdot;0)\ \hbox{ in }\R^N\backslash B(0,R).$$
Together with~\eqref{ptau1}, one concludes that $p(\tau+1,\cdot;0)\ge\sigma\,p(\tau,\cdot;0)$ in $\R^N$ and the proof of Proposition~\ref{pro1} is thereby complete.\hfill$\Box$


\SE{Proof of Theorem~\ref{th2}}\label{sec9}

The key-point in the proof of Theorem~\ref{th2} is the following result of independent interest on some monotonicity properties of the solutions of a boundary value problem in a half-line for a homogeneous linear one-dimensional reaction-diffusion.

\begin{pro} \label{propsigne}
Let $a$ and $\lambda$ be two positive real numbers and let $u$ be the solution of
\begin{equation} \label{lin1}
v_t = a\,v_{xx} + \lambda\,v,\ \ t>0,\ x>0, 
\end{equation}
with boundary condition
\begin{equation} \label{lin2}
v(t,0) = g(t),\ \ t>0,
\end{equation}
and initial datum $v_0\in L^{\infty}(0,+\infty)\backslash\{0\}$.  Assume that $v_0(x) \ge 0$ for a.e. $x>0$ and that $g$ is continuous, nonnegative and nondecreasing on $(0,+\infty)$. Then
$$
v_t(t,x) > 0
$$
provided $t\ge t_0$ and $x\ge\sqrt{8at}$, where $t_0=(2\lambda)^{-1}+e(e-1)^{-1}\lambda^{-1}>0$ is a positive constant depending only on $\lambda$.
\end{pro}

\noindent{\bf{Proof.}} Observe that $v$ can be written as
$$v=w+z,$$
where $w$ and $z$ solve~\eqref{lin1}, with $w(0,x)=v_0(x)$ and $z(0,x)=0$ for a.e. $x>0$, while $w(t,0)=0$, $z(t,0)=g(t)$ for all $t>0$.

Let us first consider the solution $w$ of the homogeneous Dirichlet boundary condition. There holds, for all $t>0$ and $x>0$,
$$w(t,x)=\int_0^{+\infty}G(t,x,y)\,v_0(y)\,dy,$$
where the Green function $G$ is given by
$$
G(t,x,y) = {e^{\lambda t} \over \sqrt{4 \pi a t} } \Bigl[ e^{ - | x - y |^2 / (4at) } - e^{- | x + y |^2 / (4at)} \Bigr]
$$
for $t>0$, $x>0$ and $y>0$. The time derivative of $G$ satisfies
$$
G_t(t,x,y) = {e^{\lambda t} \over \sqrt{4 \pi a t^3} } e^{ - | x - y |^2 / (4at) }  \Bigl[ \lambda t - {1 \over 2 }  + { | x - y |^2 \over 4 a t } \Bigr]- {e^{\lambda t} \over \sqrt{4 \pi a t^3} } e^{ - | x + y |^2 / (4at) }  \Bigl[ \lambda t - {1 \over 2}  + { | x + y |^2 \over 4 a t } \Bigr].
$$
Let us also introduce the function $\phi:[0,+\infty)\to[0,+\infty)$ defined by
$$
\phi(x) = x\,e^{-x}.
$$
This function $\phi$ has a maximum at $x=1$ and it is increasing in $[0,1]$ and decreasing in $[1,+\infty)$.

We are going to prove that $G_t(t,x,y)>0$ for all $y\in(0,+\infty)$ provided $x\ge\sqrt{8at}$ and $t\ge t_0$, where $t_0>0$ is given in the statement of Proposition~\ref{propsigne2}. In this paragraph, we fix any $t\ge t_0$, $x\ge\sqrt{8at}$ and $y>0$, and we first observe that $|x + y|^2 / (4at) \ge 2$. Two cases can then appear, depending on whether $|x-y|^2/(4at)$ is larger or smaller than $1$. In the first case, namely if $|x - y |^2 /(4at)\ge 1$, then $1 \le |x - y |^2 / (4at) < |x + y |^2 / (4at)$ and
$$
\phi \Bigl( {|x - y |^2 \over 4at } \Bigr) > \phi \Bigl( {|x + y |^2 \over 4at } \Bigr),
$$
hence we get $G_t(t,x,y) > 0$ (we also use the fact that $t \ge t_0\ge(2 \lambda)^{-1}$). In the second case, one has $0\le|x - y |^2 / (4at) <1$ and
$$
\big( \lambda t - {1 \over 2 }\big)  \bigl[ e^{ - | x - y |^2 / (4at) }  -  e^{ - | x + y |^2 / (4at) } \bigr]
\ge 
\big( \lambda t - {1 \over 2 } \big)  ( e^{-1} - e^{- 2} ) \ge e^{-1}
$$
since $t\ge t_0=(2\lambda)^{-1}+e(e-1)^{-1}\lambda^{-1}\ge(2\lambda)^{-1}$, whereas
$$
{ | x + y |^2 \over 4at}e^{ - | x + y |^2 / (4at) } < \phi(1)=e^{-1}.
$$
Hence, we also have $G_t(t,x,y) >0$ in the second case. As a consequence, $G_t(t,x,\cdot)>0$ in $(0,+\infty)$ for all $t\ge t_0$ and $x\ge\sqrt{8at}$, hence $w_t(t,x)>0$ for all $t\ge t_0$ and $x\ge\sqrt{8at}$, since $v_0$ is nonnegative and nontrivial in $(0,+\infty)$.

Let us now turn to the solution $z$ of equation~(\ref{lin1}) with vanishing initial condition and with boundary condition~(\ref{lin2}). Since $g$ is nonnegative and bounded in any interval $(0,T)$ with $T\in(0,+\infty)$, it follows from the maximum principle that $z$ is nonnegative and bounded in $(0,T)\times(0,+\infty)$ too. Furthermore, for any $h>0$ and $t>0$, there holds $z(h,\cdot)\ge 0=z(0,\cdot)$ in $(0,+\infty)$ and $z(t+h,0)=g(t+h)\ge g(t)=z(t,0)$. Hence, the maximum principle yields $z(t+h,x)\ge z(t,x)$ for all $t>0$ and $x>0$. In other words, the function $z$ is nondecreasing with respect to $t$.

As a conclusion, the solution $v$ of~(\ref{lin1}) with boundary condition~(\ref{lin2}) satisfies $v_t(t,x)> 0$ provided $t\ge t_0$ and $x\ge\sqrt{8at}$.\hfill$\Box$\break

From the previous proposition and a change of variable $x\to-x$, the following result immediately follows.

\begin{pro} \label{propsigne2}
Let $a$ and $\lambda$ be two positive real numbers and let $u$ be the solution of
\begin{equation} \label{lin1b}
v_t = a\,v_{xx} + \lambda\,v,\ \ t>0,\ x<0, 
\end{equation}
with boundary condition~\eqref{lin2} and initial datum $v_0\in L^{\infty}(-\infty,0)\backslash\{0\}$.  Assume that $v_0(x) \ge 0$ for a.e. $x<0$ and that $g$ is continuous, nonnegative and nondecreasing on $(0,+\infty)$. Then $v_t(t,x) > 0$ provided $t\ge t_0$ and $x\le-\sqrt{8at}$, where $t_0=(2\lambda)^{-1}+e(e-1)^{-1}\lambda^{-1}>0$.
\end{pro}

With these two propositions in hand, let us now turn to the proof of Theorem~\ref{th2}.\hfill\break

\noindent{\bf{Proof of Theorem~\ref{th2}.}} Let $A$, $f$, $f^{\pm}$, $\lambda^{\pm}>0$, $\theta\in(0,1)$ and $u_0$ be as in the statement and let $R>0$ and $a^{\pm}\in(0,+\infty)$ be such that
\be\label{defR2}
f(x,\cdot)=f^{\pm}\hbox{ in }[0,1]\ \hbox{ and }\ A(x)=a^{\pm}\ \hbox{ for all }|x|\ge R\hbox{ with }\pm x>0.
\ee
Let us call
\be\label{defT0}
T_0=\max\big((2\lambda^-)^{-1}+e(e-1)^{-1}(\lambda^-)^{-1},(2\lambda^+)^{-1}+e(e-1)^{-1}(\lambda^+)^{-1}\big)>0.
\ee\par
Firstly, we claim that there exists $\eta\in(0,\theta)$ such that
\be\label{claimeta}
\forall\,t\ge 1,\ \forall\,x\in\R,\ \ \big(u(t,x)\le\eta\big)\ \Longrightarrow\ \big(u(s,x)\le\theta\hbox{ for all }s\in[t,t+T_0]\big).
\ee
Assume not. Then there are a sequence $(t_n)_{n\in\N}$ in $[1,+\infty)$ and a sequence $(x_n)_{n\in\N}$ in $\R$ such that
$$u(t_n,x_n)\to0\hbox{ as }n\to+\infty\ \hbox{ and }\ \max_{[t_n,t_n+T_0]}u(\cdot,x_n)>\theta\hbox{ for all }n\in\N.$$
If the sequence $(x_n)_{n\in\N}$ were bounded, then the sequence $(t_n)_{n\in\N}$ would be bounded too due to~\eqref{defc}. Hence, up to extraction of a subsequence, $(t_n,x_n)$ would converge to a point $(t,x)$ in $[1,+\infty)\times\R$ with $u(t,x)=0$, which is impossible due to~\eqref{0u1}. Therefore, the sequence $(x_n)_{n\in\N}$ is unbounded. Up to extraction of a subsequence, let us then assume without loss of generality that $x_n\to+\infty$ as $n\to+\infty$ (the case $\lim_{n\to+\infty}x_n=-\infty$, up to extraction of a subsequence, can be handled similarly). From standard parabolic estimates, up to extraction of a subsequence, the functions $u_n:(t,x)\mapsto u(t+t_n,x+x_n)$ converge in $C^{1,2}_{t,x}$ locally in (at least) $(-1,+\infty)\times\R$ to a solution $v$ of
$$v_t=a^+v_{xx}+f^+(v)\ \hbox{ in }(-1,+\infty)\times\R$$
such that $v(0,0)=0$ and $\max_{[0,T_0]}v(\cdot,0)\ge\theta$. Furthermore, $0\le v\le1$ in $(-1,+\infty)\times\R$ by~\eqref{0u1}. Hence, $v=0$ in $(-1,0]\times\R$ from the strong maximum principle, and $v=0$ in $[0,+\infty)\times\R$ by the uniqueness of the solutions of the associated Cauchy problem. This contradicts the property $\max_{[0,T_0]}v(\cdot,0)\ge\theta\ (>0)$. Finally, the claim~\eqref{claimeta} has been proved.\par
Secondly, since the function $u$ is positive in $(0,+\infty)\times\R$ and the function $f$ is Lipschitz continuous with respect to $u\in[0,1]$ uniformly in $x\in\R$, it follows from Harnack inequality that there is a constant $C\in(0,1)$ such that, for all $(t,x)\in[1,+\infty)\times\R$,
\be\label{harnack}
u(t+T_0,x\pm\sqrt{8a^{\pm}T_0})\ge C\,u(t,x).
\ee
Let us denote
\be\label{defeps3}
\epsilon=C\,\eta
\ee
and notice that $0<\epsilon<\eta<\theta<1$. From Theorem~\ref{th1}, there is $T_{\epsilon}>0$ such that~\eqref{Teps} holds, that is,
\be\label{Teps2}
\forall\,(t,x)\in[T_{\epsilon},+\infty)\times\R,\ \ u(t,x)\ge\epsilon\ \Longrightarrow\ u_t(t,x)>0.
\ee
From~\eqref{defc}, since $\eta<\theta<1$, one can also assume without loss of generality that $T_{\epsilon}\ge1$ and that
$$\min_{[-R,R]}u(T_{\epsilon},\cdot)\ge\eta,$$
where $R>0$ is given in~\eqref{defR2}. Since $\eta>0$ and $u(T_{\epsilon},\pm\infty)=0$ by~\eqref{conv0}, it follows from the continuity of $u(T_{\epsilon},\cdot)$ that there are some real numbers $x^{\pm}$ such that
$$x^-\le-R<R\le x^+,\ \ u(T_{\epsilon},x^{\pm})=\eta\ \hbox{ and }\ u(T_{\epsilon},\cdot)\le\eta\hbox{ in }(-\infty,x^-]\cup[x^+,+\infty).$$\par
Let now $v$ be the solution of~\eqref{lin1} with $a=a^+$, $\lambda=\lambda^+$ and initial and boundary conditions given by
$$v_0=u(T_{\epsilon},\cdot+x^+)\,(>0)\hbox{ in }(0,+\infty)\ \hbox{ and }\ v(t,0)=u(t+T_{\epsilon},x^+)\,(>0)\hbox{ for all }t>0.$$
Since $u(T_{\epsilon},x^+)=\eta>\epsilon$, one infers from~\eqref{Teps2} that the continuous nonnegative function $t\mapsto u(t+T_{\epsilon},x^+)$ is actually increasing in $[0,+\infty)$. Therefore, as $T_0\ge(2\lambda^+)^{-1}+e(e-1)^{-1}(\lambda^+)^{-1}$ by~\eqref{defT0}, it follows from Proposition~\ref{propsigne} that, in particular,
$$v_t(T_0,x)>0\hbox{ for all }x\ge\sqrt{8a^+T_0}.$$
On the other hand, since $T_{\epsilon}\ge1$ and $u(T_{\epsilon},\cdot)\le\eta$ in $[x^+,+\infty)$, it follows from~\eqref{claimeta} that $u\le\theta$ in $[T_{\epsilon},T_{\epsilon}+T_0]\times[x^+,+\infty)$. In that set, since $x^+\ge R$, there holds $f(u)=\lambda^+u$ (and $A(x)=a^+$), by~\eqref{defR2} and the assumption on $f^+$. Therefore, $u(\cdot+T_{\epsilon},\cdot+x^+)$ satisfies the same linear equation as $v$ in the set $[0,T_0]\times[0,+\infty)$, with the same initial and boundary conditions on $\{0\}\times[0,+\infty)$ and $[0,T_0]\times\{0\}$. Thus,
$$u(t+T_{\epsilon},x+x^+)=v(t,x)\hbox{ for all }(t,x)\in[0,T_0]\times[0,+\infty)$$
and
\be\label{monotone}
u_t(T_0+T_{\epsilon},x)>0\hbox{ for all }x\ge x^++\sqrt{8a^+T_0}.
\ee\par
Furthermore, since $T_{\epsilon}\ge1$ and $u(T_{\epsilon},x^+)=\eta$, it follows from~\eqref{harnack} and~\eqref{defeps3} that
$$u(T_{\epsilon}+T_0,x^++\sqrt{8a^+T_0})\ge C\,\eta=\epsilon.$$
Hence, $u_t(t,x^++\sqrt{8a^+T_0})>0$ for all $t\ge T_{\epsilon}+T_0$ by~\eqref{Teps2}. Together with~\eqref{monotone} and the fact that the equation~\eqref{eq} does not depend on time, one concludes from the maximum principle applied to $u_t$ that
\be\label{monotone2}
u_t(t,x)>0\hbox{ for all }(t,x)\in[T_{\epsilon}+T_0,+\infty)\times[x^++\sqrt{8a^+T_0},+\infty).
\ee
Similarly, by using Proposition~\ref{propsigne2} among other things, one can show that
\be\label{monotone3}
u_t(t,x)>0\hbox{ for all }(t,x)\in[T_{\epsilon}+T_0,+\infty)\times(-\infty,x^--\sqrt{8a^-T_0}].
\ee\par
Finally, due to~\eqref{defc}, there is $\tau\ge T_{\epsilon}+T_0$ such that $u(\tau,\cdot)\ge\epsilon$ in $[x^--\sqrt{8a^-T_0},x^++\sqrt{8a^+T_0}]$. Hence, by~\eqref{Teps2}, there holds
$$u_t(t,x)>0\hbox{ for all }(t,x)\in[\tau,+\infty)\times[x^--\sqrt{8a^-T_0},x^++\sqrt{8a^+T_0}].$$
Together with~\eqref{monotone2} and~\eqref{monotone3}, one concludes that
$$u_t(t,x)>0\hbox{ for all }(t,x)\in[\tau,+\infty)\times\R$$
and the proof of Theorem~\ref{th2} is thereby complete.\hfill$\Box$



\begin{thebibliography}{AAA}
\footnotesize{
\bibitem{a} D.G. Aronson, {\it Bounds on the fundamental solution of a parabolic equation}, Bull. Amer. Math. Soc. {\bf 73} (1967), 890-896.
\bibitem{aw} D.G. Aronson, H.F. Weinberger, {\it Multidimensional nonlinear diffusions arising in population genetics}, Adv. Math. {\bf 30} (1978), 33-76.
\bibitem{bh} H.~Berestycki, F. Hamel, {\it Generalized transition waves and their properties}, Comm. Pure Appl. Math. {\bf 65} (2012), 592-648.
\bibitem{bhn} H. Berestycki, F. Hamel, G. Nadin, {\it Asymptotic spreading in heterogeneous diffusive media}, J.~Funct. Anal. {\bf 255} (2008), 2146-2189.
\bibitem{bhr} H.~Berestycki, F. Hamel, L.~Rossi, {\it Liouville type results for semilinear elliptic equations in unbounded domains}, Ann. Mat. Pura Appl. {\bf 186} (2007), 469-507.
\bibitem{d} E.B. Davies, {\it Heat Kernels and Spectral Theory}, Cambridge Univ. Press, 1989.
\bibitem{dm} Y. Du, H. Matano, {\it Convergence and sharp thresholds for propagation in nonlinear diffusion problems}, J.~Europ. Math. Soc. {\bf 12} (2010), 279-312.
\bibitem{fs} E.B. Fabes, D.W. Stroock, {\it A new proof of Moser's parabolic Harnack inequality using the old ideas of Nash}, Arch. Ration. Mech. Anal. {\bf 96} (1986), 327-338.
\bibitem{f} R.A. Fisher, {\it The advance of advantageous genes}, Ann. Eugenics {\bf 7} (1937), 335-369.
\bibitem{fr} A. Friedman, {\it Partial Differential Equations of Parabolic Type}, Prentice-Hall, Englewood Cliffs, New Jersey, 1964.
\bibitem{ggn} J. Garnier, T. Giletti, G. Nadin, {\it Maximal and minimal spreading speeds for reaction-diffusion equations in nonperiodic slowly varying media}, J.~Dyn. Diff. Equations {\bf 24} (2012), 521-538.
\bibitem{hn} F. Hamel, G. Nadin, {\it Spreading properties and complex dynamics for monostable reaction-diffusion equations}, Comm. Part. Diff. Equations {\bf 37} (2012),  511-537.
\bibitem{kpp} A.N. Kolmogorov, I.G. Petrovsky, N.S. Piskunov, {\it \'Etude de l'\'equation de la diffusion avec croissance de la quantit\'e de mati\`ere et son application \`a un probl\`eme biologique}, Bull. Univ. \'Etat Moscou, S\'er. Intern.~A {\bf 1} (1937), 1-26.
\bibitem{m} J.D. Murray, {\it Mathematical Biology}, Springer-Verlag, 2003.
\bibitem{n} J. Norris, {\it Long-time behaviour of heat flow: global estimates and exact asymptotics}, Arch. Ration. Mech. Anal. {\bf 140} (1997), 161-195.
\bibitem{r} J.-M. Roquejoffre, {\it Eventual monotonicity and convergence to travelling fronts for the solutions of parabolic equations in cylinders}, Ann. Inst. H. Poincar\'e, Anal. Non Lin\'eaire {\bf 14} (1997), 499-552.
\bibitem{w} D.V. Widder, {\it The Heat Equation}, Academic Press, 1975.
\bibitem{y} E. Yanagida, {\it Irregular behavior of solutions for Fisher's equation}, J.~Dyn. Diff. Equations {\bf 19} (2007), 895-914.
\bibitem{z} A. Zlato{\v{s}}, {\it Propagation of reactions in inhomogeneous media}, preprint, 2014.
}
\end{thebibliography}
\end{document}